\documentclass[11pt]{amsart}
\usepackage{geometry}                % See geometry.pdf to learn the layout options. There are lots.
\usepackage{graphicx}
\usepackage{amssymb}
\usepackage{epstopdf}
\usepackage{subfig} % extra packages for figures and Taylor diagrams
\usepackage{color}
\usepackage{caption}
\usepackage{cancel}
\newcommand {\R} {\mathbb {R}}

\newcommand {\Z} {\mathbb {Z}}

\newcommand {\Pj} {\mathbb {P}}

\newcommand {\z} {\textbf{0}}
\newcommand {\no} {\noindent}

\newcommand {\pa} {\partial}
\newcommand {\rest}[1] {\big|_{#1}}
\theoremstyle{definition}

\newtheorem{lemma}{Lemma}
\newtheorem{defn}{Definition}[section]

\newtheorem{example}{Example}
\newtheorem{thm}{Theorem}[section]

\newtheorem{eg}[thm]{Example}
\newtheorem{rem}[thm]{Remark}

\title{The title of my page}

\title{A Monster Tower Approach to Goursat Multi-flags}
\author{Alex L. Castro}
\email{alex.castro@utoronto.ca}
\address{University of Toronto}
\author{Wyatt C. Howard}
\email{whoward@ucsc.edu}
\thanks{W.C.H. is the corresponding author}
\address{University of California at Santa Cruz}
\date{\today}     
\begin{document}

% ----------------------------------------------------------------

\maketitle

\begin{abstract}

We consider here the problem of classifying orbits of an action of the diffeomorphism group of 3-space on a tower of fibrations with $\mathbb{P}^2$-fibers that generalize the Monster Tower due to Montgomery and Zhitomirskii. As a corollary we give the first steps towards the problem of classifying Goursat 2-flags of small length.  In short, we classify the orbits within the first four levels of the Monster Tower and show that
there is a total of $34$ orbits at the fourth level in the tower.     

% An abstract is a short summary of the main points in the project.  It
% does not mention background material or detail.

\end{abstract}

% ----------------------------------------------------------------
%\tableofcontents
%\addtocontents{toc}{\protect\setcounter{tocdepth}{1}}
% ----------------------------------------------------------------

%%%%%%%%%%%%%%%%%%%%%%%%%%%%%%%%%%%%%%%%%%%%%%%%%%%%%%%%%%%%%%%%%%%%%%%%%%%%%%%%%%%%%%%%%%%%%%%%%%%%%%%%%%%%%%%

\section{Introduction}\label{sec:intro}

A Goursat flag is a nonholonomic distribution $D$ with ``slow growth''.
By slow growth we mean that the associated flag of distributions 
$$D \hspace{.2in} \subset  \hspace{.2in} D + [D,D]\hspace{.2in} \subset  \hspace{.2in} D + [D,D] + [[D,D],[D,D]] \dots ,$$
grows by one dimension at each bracketing step and after $n$ steps it will span the entire tangent bundle.  By an abuse of notation, $D$ in this context also means the sheaf of vector fields spanning $D$.

Though less popular than her nonholonomic siblings like the contact distribution, or rolling distribution in mechanics \cite{gil}, 
Goursat distributions are more common than one would think.  The canonical Cartan distributions in the jet spaces 
$J^k(\mathbb{R},\mathbb{R})$, or the non-slip constraint for a jackknifed truck \cite{jean} are examples. 

Generalizations of Goursat flags have been proposed in the literature. One such notion is that of a {\it Goursat multi-flag.} Typical examples of Goursat multi-flags include the Cartan distributions $C$ of the jet spaces $J^k(\mathbb{R},\mathbb{R}^n),n\geq 2.$ 
Iterated bracketing this time:
$$C\hspace{.2in} \subset  \hspace{.2in} C + [C,C] \hspace{.2in} \subset  \hspace{.2in} C + [C,C] + [[C,C],[C,C]] \dots,$$ 
leads to a jump in rank by $n$ dimensions at each step. 

To our knowledge the general theory behind Goursat multi-flags made their first appearance in the works of A. Kumpera and J. L Rubin \cite{kumpera1}. 
P. Mormul has also been very active in breaking new ground \cite{mormul1}, and developed new combinatorial tools to investigate the normal forms of these distributions.  In addition to Mormul's work is Yamaguchi's work \cite{yamaguchi1} and \cite{yamaguchi2} where he investigated the local properties of Goursat multi-flags.  It is also important to mention that results from his work in \cite{yamaguchi2} were crucial
to our classification proceedure.      

In this note we concentrate on the problem of classifying local germs of Goursat multi-flags of small length.  
We will consider Goursat 2-flags of length up to 4.  Goursat 2-flags exhibit many new geometric features our old Goursat 
flags(Goursat $1$-flags) did not possess. The geometric properties of Goursat multi-flags was the main subject of the paper \cite{castro}.

We approach the classification problem from a geometric standpoint, and will partially follow the programme started by Montgomery and Zhitomirskii in \cite{mont1}. 
Our starting point is the universality theorem proved by Yamaguchi and Shibuya \cite{yamaguchi1} stating that any germ $D$ of a Goursat distribution, or using Mormul's terminology, a Goursat $n$-flag of length $k$ is analytically equivalent to the germ of the canonical distribution on 
the $k$-step Cartan prolongation \cite{bryant} of the flat $\mathbb{R}^{n+1}$ with its trivial bundle, the germs being taken at some
$D$-dependent point. 

In what follows we apply Yamaguchi's theorem to the $n=2$ case, and translate the classification problem of Goursat 2-flags into a problem of classifying points in a tower of real $\mathbb{P}^2$ (projective plane) fibrations 

\begin{equation}\label{eqn:tower}
\cdots \rightarrow \mathcal{P}^{4}(2)\rightarrow \mathcal{P}^{3}(2)\rightarrow \mathcal{P}^{2}(2) \rightarrow 
\mathcal{P}^{1}(2) \rightarrow \mathcal{P}^{0}(2) = \mathbb{R}^3,
\end{equation}
where $\text{dim}(\mathcal{P}^{k}(2)) = 3 + 2k$ and $\mathcal{P}^{k}(2)$ is the Cartan prolongation of $\mathcal{P}^{k-1}(2)$.  The global topology of these manifolds is much more interesting though, and has yet to be explored (\cite{castro}). 

Each $\mathcal{P}^{k}(2)$ comes equipped with a rank-3 distribution $\Delta_k$.  At a dense open subset of 
$\mathcal{P}^{k}(2)$ this $\Delta_{k}$ is locally diffeomorphic to the Cartan distribution $C$ in $J^k(\mathbb{R},\mathbb{R}^2).$ 
$\Delta_k$ has an associated flag of length $k$. 
  
A description of symmetries of the $\Delta_k$ is the content of a theorem due to Yamaguchi \cite{yamaguchi2}, attributed 
to Backl\"und, and dating from the $1980$'s. He showed that that any symmetry of $(\mathcal{P}^{k}(2),\Delta_k)$ results from the Cartan prolongation of a symmetry of the base manifold, i.e. of a diffeomorphism of the 3-fold $\mathbb{R}^3$. 

By applying the techniques developed in \cite{castro,mont2}, we will attack the classification problem utilizing the {\it curves-to-points} approach and a new technique we named {\it isotropy representation}, and used to some extent in \cite{mont1}, somehow inspired by \'E. Cartan's moving frame method \cite{favard}. Our main result states that there are $34$ Goursat $2$-flags of length 4, and we provide the exact numbers for each length size. Our approach is constructive.  Normal forms for each equivalence class can be made explicit. Due to space limitations we will write down only a couple instructive examples.

We would like to mention that P. Mormul and Pelletier \cite{mormul2} have attempted an alternative solution to the classification problem.
In their classification work, they employed the results and tools proved from previous works done by Mormul.  In \cite{mormul} Mormul discusses two coding systems for special
$2$-flags and proves that the two coding systems are the same.  One system is the Extend Kumpera Ruiz system, which is a coding of system
used to describe $2$-flags.  The other is called Singularity Class coding, which is an intrinsic coding system that describes
the sandwich diagram \cite{mont1} associated to $2$-flags.  A brief outline on how these coding systems relate
to Montgomery's $RVT$ coding is discussed in \cite{castro}.  Then, building upon Mormul's work in \cite{mormul3}, Mormul and Pelletier 
use the idea of strong nilpotency of special multi-flags, along with properties and relationships between his two coding systems, to classify the
these distributions, up to length $4$(equivalently, up to level $4$ of the Monster Tower).  Our $34$ agrees with theirs.

Here is a short description of the paper. In section one we acquaint ourselves with the main definitions necessary for the statements of our main results. Section two contains our precise statements, and a few explanatory remarks to help the reader progress through the theory with with us.
Section three consists of the statements of our main results.  In section four we discuss the basic tools and ideas that will
be needed to prove our various results. Section five is devoted to technicalities, and the actual proofs. We conclude the paper, section six, with a quick summary of our findings. 

For the record, we have also included an appendix where our lengthy computations are scrutinized.  

\vspace{.2in}

\noindent {\bf Acknowlegements.} We would like to thank Corey Shanbrom, and Richard Montgomery (both at UCSC) for many useful conversations and remarks.

%-----------------------------------------------------------------
\section{Main definitions}

\subsection{Constructions} 
Let $(Z,\Delta)$ be a manifold of dimension $n$ equipped with a plane field of rank $r$ and let $\Pj (\Delta)$ be the {\it projectivization} of $\Delta$.
As a manifold, $$Z^1 = \Pj (\Delta) .$$ Various objects in $Z$ can be canonically prolonged  (lifted) to the new manifold $Z^1$. 

\begin{table}\label{tab:prol}
\caption{Some geometric objects and their Cartan prolongations.}
\begin{tabular}{|c|c|}
  \hline
   curve $c:(\mathbb{R},0)\rightarrow (Z,p)$ & curve $c^{1}:(\mathbb{R},0)\rightarrow (Z^1,p),$ \\
   & $c^{1}(t) =\text{(point,moving line)} = (c(t),span\{ \frac{dc(t)}{dt} \}) $\\ \hline
  diffeomorphism $\Phi: Z \circlearrowleft$ & diffeomorphism $\Phi^{1}: Z^1 \circlearrowleft$, \\
  & $\Phi^{1}(p,\ell) = (\Phi(p),d\Phi_p(\ell))$ \\ \hline
  rank $r$ linear subbundle  & rank $r$ linear subbundle $\Delta_1=d\pi_{(p,\ell)}^{-1}(\ell)\subset TZ^1$,\\
  $\Delta \subset TZ$ & $\pi: Z^1 \rightarrow Z$ is the canonical projection. \\
  \hline
\end{tabular}
\end{table}

Given an analytic curve $c:(I,0)\rightarrow (Z,p)$, where $I$ is some open interval about zero in $\mathbb{R}$, we can naturally define a curve a new curve $$c^{1}:(I,0)\rightarrow (Z^1,(p,\ell))$$ with image in $Z^1$ and where $\ell = span\{ \frac{dc(0)}{dt} \}$. 
This new curve, $c^{1}(t)$, is called the \textit{prologation} of $c(t)$.  If $t = t_{0}$ is not a regular point, then we define $c^{1}(t_{0})$ by taking the limit $\lim_{t \rightarrow t_{0}} c^{1}(t)$ where the limit varies over the regular points $t \rightarrow t_{0}$.  An important fact to
note, proved in \cite{mont1}, is that the analyticity of $Z$ and $c$ implies that the limit is well defined and that the prolonged curve
$c^{1}(t)$ is analytic as well.  Since this process can be iterated, we will write 
$c^{k}(t)$ to denoted the $k$-fold prolongation of the curve $c(t)$.
  
The manifold $Z^1$ also comes equipped with a distribution $\Delta_1$ called the {\it Cartan prolongation of $\Delta$} \cite{bryant} and is defined as follows. 
Let $\pi : Z^1 \rightarrow Z$ be the projection map $(p, \ell)\mapsto p$. Then
$$\Delta_1(p,\ell) = d\pi_{(p,\ell)}^{-1}(\ell),$$
i.e. {\it it is the subspace of $T_{(p,\ell)}Z^1$ consisting of all tangents to curves which are prolongations of curves one level down through $p$ in the direction $\ell$.}  
It is easy to check using linear algebra that $\Delta_1$ is also a $k$-plane field. The pair $(Z^1,\Delta_{1})$ is called the {\it Cartan prolongation} of $(Z,\Delta)$. 

\begin{example}
Take $Z = \mathbb{R}^{3}$ with its tangent bundle, which we denote by $\Delta_0$. 
Then the tower shown in equation (\ref{eqn:tower}) is obtained by prolonging the pair $(\mathbb{R}^{3},\Delta_0)$ four times. 
\end{example}
%Note that locally, $P_n^k \approx \mathbb{R}^n \times (\mathbb{P}^{n-1})^k$ and $\dim (P_n^k) = n + k(n-1)$.  
By a {\it symmetry} of the pair $(Z,\Delta)$ we mean a (local or global) diffeomorphism $\Phi \in Diff(3)$ of $Z$ that preserves the subbundle $\Delta$.
The symmetries of $(Z,\Delta)$ can also be prolonged to symmetries $\Phi^{1}$ of $(Z^1,\Delta_{1})$ as follows. 
Define 
$$\Phi^{1}(p,\ell)=(\Phi(p), d\Phi(\ell)).$$ Since $d\Phi(p)$ is invertible, and $d\Phi$ is linear the second component is well defined as a projective map. The new born symmetry we denoted by $\Phi^{1}$ is the prolongation (to $Z^1$) of $\Phi$. 
Objects of interest in this paper and their Cartan prolongations are summarized in table (1). Unless otherwise mentioned prolongation will always refer to Cartan prolongation.

\begin{eg}[Prolongation of a cusp.]
Let $c(t) = (t^{2}, t^{3}, 0)$ be the $A_{2}$ cusp in $\R^{3}$.  Then
$c^{1}(t) = (x(t), y(t), [dx, dy, dz]) = (t^{2}, t^{3}, 0, [2t: 3t^{2}: 0])$.  After we introduce fiber affine coordinates
$u = \frac{dy}{dx}$ and $v = \frac{dz}{dx}$ around the point $[1: 0: 0]$ we obtain the immersed curve
$$c^{1}(t) = (t^{2}, t^{3}, 0, \frac{3}{2}t, 0)$$
\end{eg}

\subsection{Constructing the Monster tower.}
We start with $\R^{n+1}$ as our base and take $\Delta_{0} = T \R^{n+1}$.  Prolonging $\Delta_{0}$ we get that 
$\mathcal{P}^{1}(n) = \Pj \Delta_{0}$  with the distribution $\Delta_{1}$.  By
iterating this process we end up with the manifold $\mathcal{P}^{k}(n)$ which is endowed with the rank $n$
distribtuion $\Delta_{k} = (\Delta_{k-1})^{1}$ and fibered over $\mathcal{P}^{k-1}(n)$.
In this paper we will be looking at the case of when $n=3$.

\begin{defn}
The Monster tower is a sequence of manifolds with distributions, $(\mathcal{P}^{k}, \Delta_{k})$, 
together with fibrations $$\cdots \rightarrow \mathcal{P}^{k}(n) \rightarrow \mathcal{P}^{k-1}(n) \rightarrow \cdots \rightarrow 
\mathcal{P}^{1}(n) \rightarrow \mathcal{P}^{0}(n) = \mathbb{R}^{n+1}$$
and we write $\pi_{k,i}: \mathcal{P}^{k}(n) \rightarrow \mathcal{P}^{i}(n)$, with $i < k$ for the projections.
\end{defn}

\begin{thm}
For $n > 1$ and $k>0$ any local diffeomorphism of $\mathcal{P}^{k}(n)$ preserving the distribution $\Delta_{k}$
is the restriction of the $k$-th prolongation of a local diffeomorphism $\Phi \in Diff(n)$.
\end{thm}

Proof: This was shown by Yamaguchi and Shibuya in (\cite{yamaguchi2}).

\begin{rem}
The importance of the above result cannot be stressed enough.  
This theorem by Yamaguchi and Shibuya is what allows the isotropy representation method, discussed in section five of the paper, 
of classifying orbits within the Monster Tower to work.
\end{rem}

\begin{rem}
Since we will be working exclusively with the $n=2$ Monster tower in this paper, we will just write
$\mathcal{P}^{k}$ for $\mathcal{P}^{k}(2)$.
\end{rem}

\begin{defn}
Two points $p,q$ in $\mathcal{P}^k$ are said to be equivalent if and only if there is a $\Phi\in \text{Diff}(3)$ such that $\Phi^{k}(p)=q$, in other words, $q\in \mathcal{O}(p)$($\mathcal{O}(p)$ is the orbit of the point $p$). 
\end{defn}

\subsection{Orbits.}
Yamaguchi's theorem states that any symmetry of $\mathcal{P}^k$ comes from prolonging a diffeomorphism of $\mathbb{R}^3$ $k$-times. This remark is essential to our computations. Let us denote by $\mathcal{O}(p)$ the orbit of the point $p$ under the action of $Diff(\mathbb{R}^3)$.
In trying to calculate the various orbits within the Monster tower we see that it is easier to fix the base points
$p_{0} = \pi_{k,0}(p_{k})$ and $q_{0} = \pi_{k,0}(q_{k})$ to be $0 \in \R^{3}$.  This means that we can replace the pseudogroup 
$Diff(3)$, diffeomorphism germs of $\R^{3}$, by the group $Diff_{0}(3)$ of diffeomorphism germs that map the origin back to the 
origin in $\R^{3}$.
    
\begin{figure}
  \centering
  \subfloat[Prolongation of the distribution.]{\label{fig:prologdist}\includegraphics[width=0.35\textwidth]{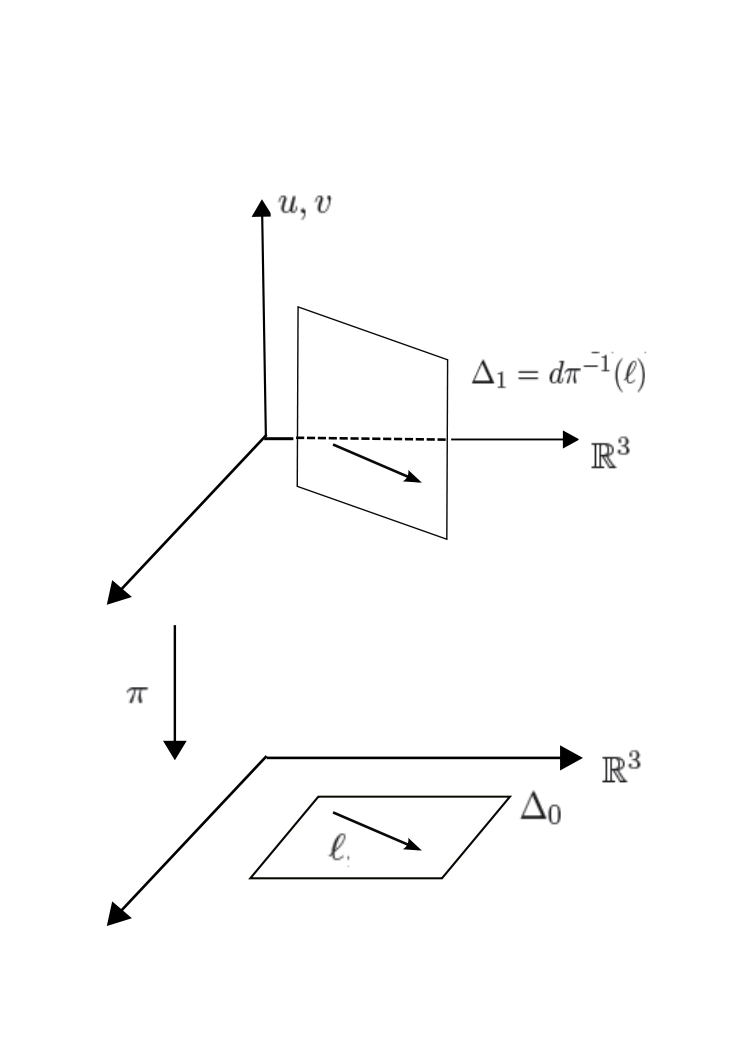}}              
  \subfloat[Prolongation of a diffeomorphism]{\label{fig:prologdiff}\includegraphics[width=0.35\textwidth]{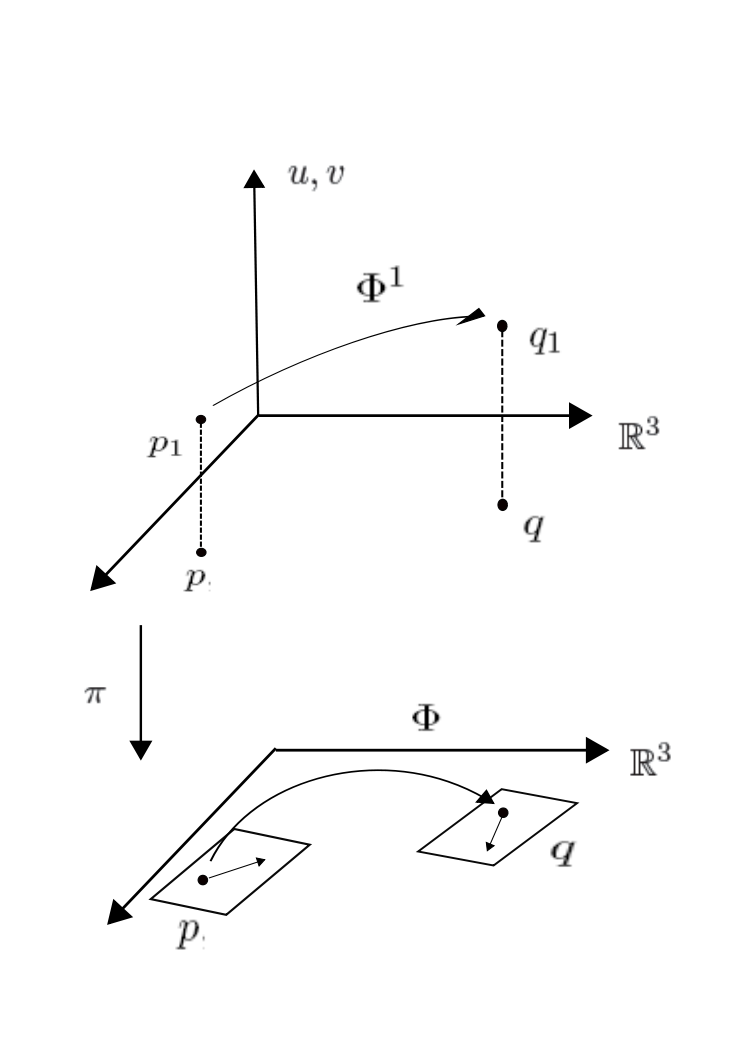}}
  \label{fig:prolongation}
\end{figure}

\begin{defn}
We say that a curve or curve germ $\gamma: (\R, 0) \rightarrow (\R^{3}, 0)$ realizes the point $p_{k} \in \mathcal{P}^{k}$
if $\gamma^{k}(0) = p_{k}$, where $p_{0} = \pi_{k,0}(p_{k})$.
\end{defn}

\begin{defn}
A direction $\ell \subset \Delta_{k}(p_{k})$, $k \geq 1$ is called a critical direction if there exists an
immersed curve, at level $k$, that is tangent to the direction $\ell$ whose projection to the zero-th level is the constant curve.
If no such curve exists,
then we call $\ell$ a regular direction.
\end{defn}

\begin{defn}
Let $p \in \mathcal{P}^{k}$, then the
\begin{eqnarray*}
\text{Germ}(p) &=& \{ c :(\mathbb{R},0)\rightarrow (\mathbb{R}^3,0)| \text{$\frac{dc^{k}}{dt}\vert_{t=0}\neq 0$ is a regular direction} \}.
\end{eqnarray*}
\end{defn}

\begin{defn}
Two curves $\gamma$, $\sigma$ in $\R^{3}$ are $(RL)$ equivalent, written $\gamma \sim \sigma$ $\Leftrightarrow$ there 
exists a diffeomorphism germ $\Phi \in Diff(3)$ and a reparametrization $\tau \in Diff_{0}(1)$ of $(\R,0)$ such that
$\sigma = \Phi \circ \gamma \circ \tau$. 
\end{defn}

% ----------------------------------------------------------------
\section{Main results}
\begin{thm}[Orbit counting per level]\label{thm:main}
In the $n=2$ Monster tower the number of orbits within each of the first four levels of the tower are as follows:
level $1$ has $1$ orbit, level $2$ has $2$ orbits, level $3$ has $7$ orbits, and level $4$ has $34$ orbits.
\end{thm}

The main idea behind this classification is a coding system developed by Castro and Montgomery \cite{castro}.  This
coding system is known as $RVT$ coding where each point in the Monster tower is labeled by a sequence of 
$R$'s, $V$'s, $T$'s, and $L$'s along with various decorations.  We will give an explanation of this 
coding system in the next section.  Using this coding system we went class by class and determined the number
or orbits within every possible $RVT$ class that could arise at each of the first four levels.

\begin{thm}[Listing of orbits within each $RVT$ code.]\label{thm:count}
The above table, is a break down of the number of orbits that appear within each $RVT$ class within the first three levels.

\begin{table}\label{tab:codes}
\begin{tabular}{|c|c|c|c|}
\hline
Level of tower &      $RVT$ code  & Number of orbits & Normal forms \\
$1$            &        $R$       & $1$              &  $(t,0,0)$ \\
$2$            &        $RR$      & $1$              &  $(t,0,0)$ \\
               &        $RV$      &  $1$             &  $(t^{2}, t^{3}, 0)$\\
$3$            &       $RRR$      & $1$              &  $(t,0,0)$ \\
               &       $RRV$      & $1$              &  $(t^{2}, t^{5}, 0)$ \\
               &       $RVR$      & $1$              &  $(t^{2}, t^{3}, 0)$ \\
               &       $RVV$      & $1$              &  $(t^{3}, t^{5}, t^{7})$, $(t^{3}, t^{5}, 0)$ \\
               &       $RVT$      & $2$              &  $(t^{3}, t^{4},t^{5})$, $(t^{3}, t^{4}, 0)$ \\
               &       $RVL$      & $1$              &  $(t^{4}, t^{6}, t^{7})$ \\
\hline
\end{tabular}
\end{table}

\no For level $4$ there is a total of $23$ possible $RVT$ classes.  Of the $23$ possibilities $14$ of them consist of a single orbit.
The classes $RRVT$, $RVRV$, $RVVR$, $RVVV$, $RVVT$, $RVTR$, $RVTV$, $RVTL$ consist of $2$ orbits, and the class $RVTT$ consists of $4$ orbits.
\end{thm}   

\begin{rem}
There are a few words that should be said to explain the normal forms column in table $2$.
Let $p_{k} \in \mathcal{P}^{k}$, for $k = 1,2,3$, having $RVT$ code $\omega$, meaning $\omega$ is a word
from the second column of the table.  For $\gamma \in Germ(p_{k})$,
then $\gamma$ is $(RL)$ equivalent to one of the curves listed in the normal forms column for the 
$RVT$ class $\omega$.  Now, notice that for the class $RVV$ that there are two inequivalent curves sitting in the 
normal forms column, but that there is only one orbit within that class.  This is because the two normal forms are equal
to each other, at $t=0$, after three prolongations.  However, after four prolongations they
represent different points at the fourth level.  This corresponds to the fact that at the fourth level class
$RVVR$ breaks up into two orbits.
\end{rem}

The following theorems are results that were proved in \cite{castro} and which helped to reduce the number calculations in our orbit classification process.

\begin{defn}
A point $p_{k} \in \mathcal{P}^{k}$ is called a Cartan point if its $RVT$ code is $R^{k}$.
\end{defn}

\begin{thm}\label{thm:cartan}
The $RVT$ class $R^{k}$ forms a single orbit at any level within the Monster tower $\mathcal{P}^{k}(n)$ for $k \geq 1$ and $n \geq 1$.
Every point at level $1$ is a Cartan point.  For $k > 1$ the set $R^{k}$ is an open dense subset of $\mathcal{P}^{k}(n)$. 
\end{thm}

\begin{defn}
A parametrized curve is an $A_{2k}$ curve, $k \geq 1$ if it is equivalent to the curve
$$(t^{2}, t^{2k + 1}, 0) $$
\end{defn}

\begin{thm}\label{thm:ak}
Let $p_{j} \in \mathcal{P}^{j}$ with $j = k + m + 1$, with $m,k \geq 0$, $m,k$ are positive integers, and $p_{j} \in R^{k}CR^{m}$,
then $Germ(p_{j})$ contains a curve germ equivalent to the $A_{2k}$ singularity, which means that the $RVT$ class
$R^{k}CR^{m}$ consists of a single orbit.
\end{thm}

\begin{rem}
One could ask ``Why curves?'' The space of $k$-jets of functions $f:\mathbb{R}\rightarrow \mathbb{R}^{2}$, usually denoted by $J^k(\mathbb{R},\mathbb{R}^{2})$ is an open dense subset of $\mathcal{P}^{k}$. It is in this sense that a point $p\in \mathcal{P}^{k}$ is roughly speaking the $k$-jet of a curve in $\mathbb{R}^3$. Sections of the bundle 
$$J^k(\mathbb{R},\mathbb{R}^{2}) \rightarrow \mathbb{R}\times \mathbb{R}^{2}$$
are $k$-jet extensions of functions. 
Explicitly, given a function $t\mapsto (t,x(t),y(t))$ its $k$-jet extension is defined as 
$$(x,f(x))\mapsto (t,x(t),y(t),x'(t),y'(t),\dots,x^{(k)}(t),y^{(k)}(t)).$$ (Superscript here denotes the order of the derivative.)
It is an instructive example to show that for certain choices of fiber affine coordinates in $\mathcal{P}^{k}$, not involving critical directions, that our local charts will look like a copy of $J^k(\mathbb{R},\mathbb{R}^{2})$.
\

Another reason for looking at curves is because it gives us a better picture for the overall behavior of an $RVT$ class.  
If one knows all the possible curve normal forms for a particular $RVT$ class, say $\omega$, then not only does one know how many
orbits are within the class $\omega$, but they also know how many orbits are within the regular prolongation of $\omega$.  By
regular prolongation of an $RVT$ class $\omega$ we mean the addition of $R$'s to the end of the word $\omega$, i.e. the regular
prolongation of $\omega$ is $\omega R \cdots R$.  This method of using curves to classify $RVT$ classes was used in \cite{mont1} and proved to be very 
successful in classifing points within the $n=1$ Monster Tower.
\end{rem}

%The Yamaguchi-Shibuya theorem then implies that,

%\begin{cor}[Goursat 2-flags by length size] The number of equivalence classes of Goursat 2-flags, of sizes $k=1,2,3,4$ are described in the table.  
%\begin{table}\label{tab:prol}
%\caption{Some geometric objects and their Cartan prolongations.}
%\begin{tabular}{|c|c|}
%  \hline
%   length of Goursat 2-flag & no. equivalence classes \\ \hline
%   1 & \\ 
%  2 & \\ 
% 3 & \\ 
%  4 & \\ 
%  \hline
%\end{tabular}
%\end{table}
%\end{cor}

\section{Tools and ideas involved in the proofs}\label{sec:tools} 
Before we begin the proofs we need to define the $RVT$ code.

\subsection{$RC$ coding of points.}

\begin{defn}
A point $p_{k} \in \mathcal{P}^{k}$, where $p_{k} = (p_{k-1}, \ell)$ is called regular or critical point if the 
line $\ell$ is a regular direction or a critical direction. 
\end{defn}

\begin{defn}
For $p_{k} \in \mathcal{P}^{k}$, $k \geq 1$ and $p_{i} = \pi_{k,i}(p_{k})$, we write
$\omega_{i}(p_{k}) = R$ if $p_{i}$ is a regular point and $\omega_{i}(p_{k}) = C$ if $p_{i}$ is a critical point.
Then the word $\omega(p_{k}) = \omega_{1}(p_{k}) \cdots \omega_{k}(p_{k})$ is called the $RC$ code for the point $p_{k}$.
Note that $\omega_{1}(p_{k})$ is always equal to $R$ by Theorem $3.4$.  
\end{defn}

\no So far we have not discussed how critical directions arise inside of $\Delta_{k}$.  The following section will show
that there is more than one kind of critical direction that can appear within the distribution $\Delta_{k}$.

\subsection{Baby Monsters.}
One can apply prolongation to any analytic $n$-dimensional manifold $F$ in place of $\R^{n}$.  Start out with 
$\mathcal{P}^{0}(F) = F$ and take $\Delta^{F}_{0} = TF$.  Then the prolongation of the pair 
$(F, \Delta^{F}_{0})$ is $\mathcal{P}^{1}(F) = \Pj TF$, with canonical rank $m$ distribution
$\Delta^{F}_{1} = (\Delta^{F}_{0})^{1}$.  By iterating this process $k$ times we end up with the pair 
$(\mathcal{P}^{k}(F), \Delta^{F}_{k})$, which is analytically diffeomorphic to $(\mathcal{P}^{k}(n-1), \Delta_{k})$.
\

Now, apply this process to the fiber $F_{i}(p_{i}) = \pi^{-1}_{i, i-1}(p_{i-1}) \subset \mathcal{P}^{i}$ through the point
$p_{i}$ at level $i$. The fiber is an $(n-1)$-dimensional integral submanifold for $\Delta_{i}$.  Prolonging, we see the
$\mathcal{P}^{1}(F_{i}(p_{i})) \subset \mathcal{P}^{i + 1}$, and has the associated distribution
$\delta^{1}_{i} = \Delta^{F_{i}(p_{i})}_{1}$; that is,
$$\delta^{1}_{i}(q) = \Delta_{i + 1}(q) \cap T_{q}(\mathcal{P}^{1}(F_{i}(p_{i}))) $$
which is a hyperplane within $\Delta_{i + 1}(q)$, for $q \in \mathcal{P}^{1}(F_{i}(p_{i}))$.  When this prolongation process is 
iterated, we end up with the submanifolds
$$\mathcal{P}^{j}(F_{i}(p_{i})) \subset \mathcal{P}^{i + j}$$
with the hyperplane subdistribution $\delta^{j}_{i} \subset \Delta_{i + j}(q)$ for $q \in \mathcal{P}^{j}(F_{i}(p_{i}))$.

\begin{defn}   A baby Monster born at level $i$ is a sub-tower $(\mathcal{P}^{i}(F_{i}(p_{i})), \delta^{j}_{i})$,
for $j \geq 0$ within the Monster tower.  If $q \in \mathcal{P}^{j}(F_{i}(p_{i}))$ then we will say that a baby Monster born
at level $i$ passes through $q$ and that $\delta^{j}_{i}(q)$ is a critical hyperplane passing through $q$, and born at level $i$. 
\end{defn}

\begin{rem} The vertical plane $V_k (q)$, which is of the form $\delta^{0}_{k} (q)$,  is always one
of the  critical hyperplanes 
passing through $q$. 
\end{rem}

\begin{thm}
A direction $\ell \subset \Delta_{k}$ is called critical $\Leftrightarrow$ $\ell$ is contained in a critical hyperplane.
\end{thm}

\begin{figure}
  \centering
  \subfloat[Above a regular point.]{\label{fig:one-plane}\includegraphics[clip = true, width=0.3\textwidth]{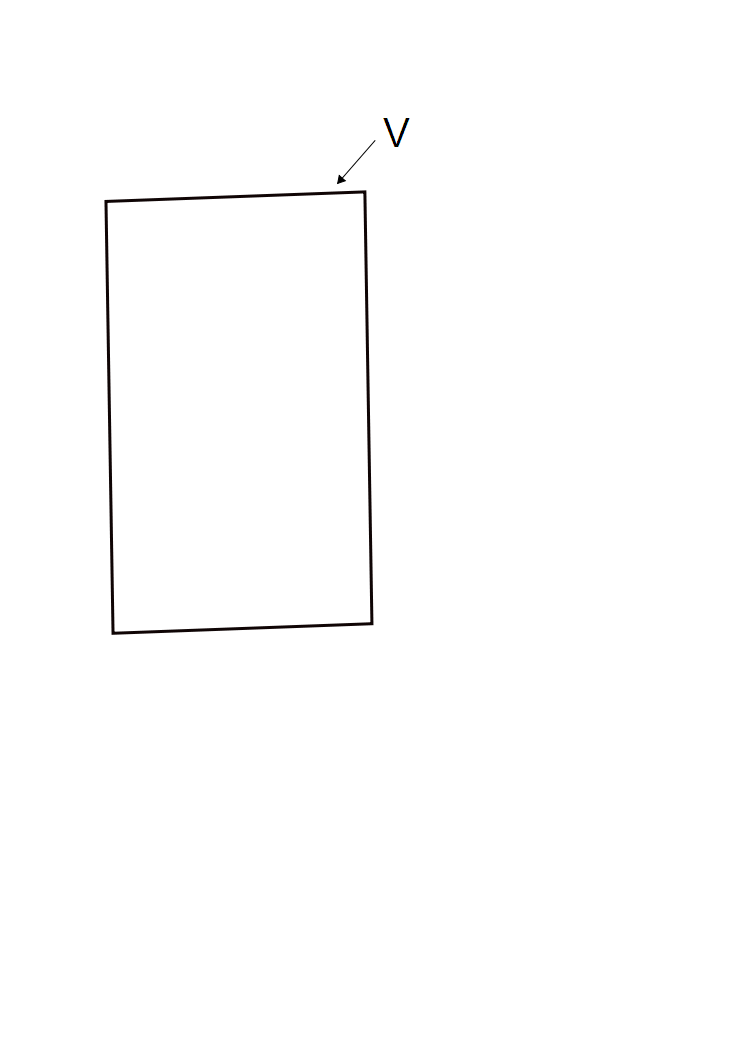}}              
  \subfloat[Above a vertical or tangency point.]{\label{fig:two-planes}\includegraphics[width=0.3\textwidth]{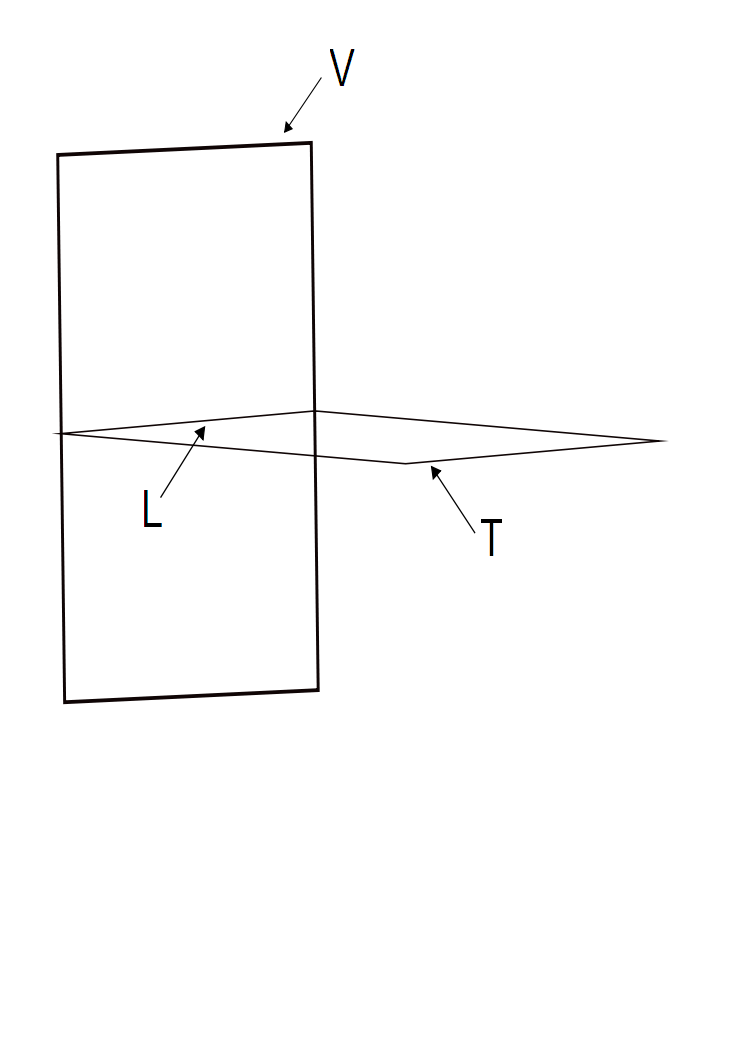}}
  \subfloat[Above an $L$ point.]{\label{fig:three-planes}\includegraphics[width=0.3\textwidth]{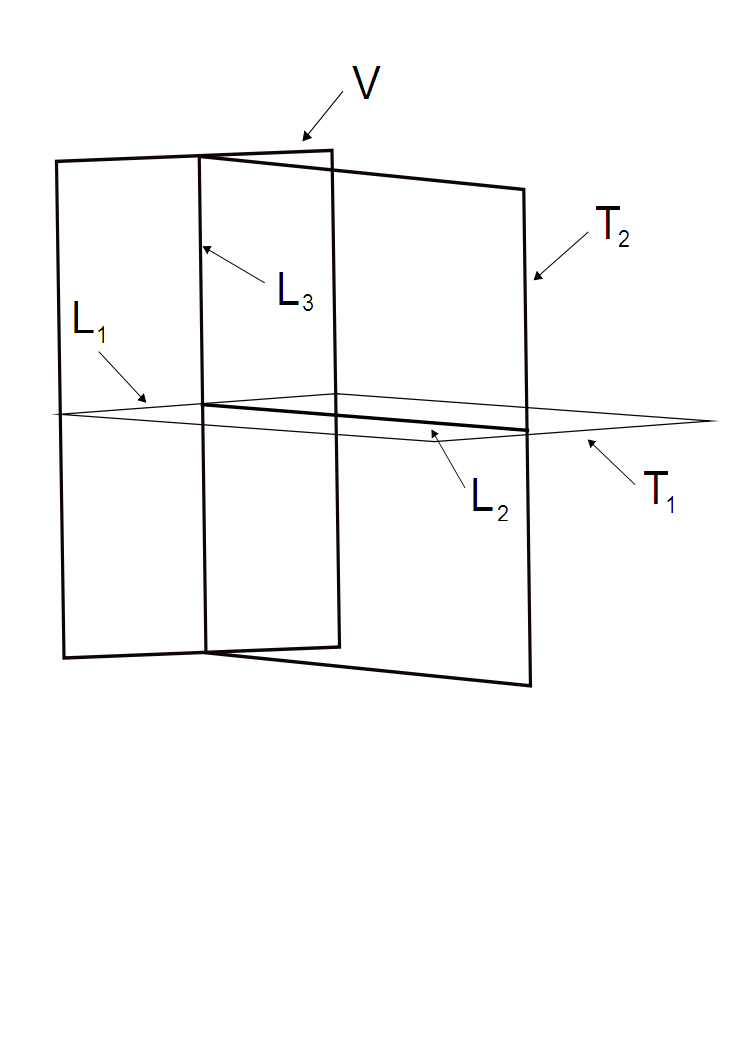}}
  \caption{Arrangement of critical hyperplanes.}
  \label{fig:arrangement}
\end{figure}

\subsection{Arrangements of critical hyperplanes for $n = 2$.}
Over any point $p_{i}$, at the i-th level of the Monster tower, there is a total of three different hyperplane configurations
for $\Delta_{i}$.  These three configurations are shown in diagrams $(a)$, $(b)$, and $(c)$.  Figure $(a)$ is the picture
for $\Delta_{i}(p_{i})$ when the i-th letter in the $RVT$ code for $p_{i}$ is the letter $R$.  From our earlier discussion, this means that the
vertical hyperplane, labeled with a $V$, is the only critical hyperplane sitting inside of $\Delta_{i}(p_{i})$.  Figure $(b)$ is the
picture for $\Delta_{i}(p_{i})$ when the i-th letter in the $RVT$ code is either the letter $V$ or the letter $T$.  This gives that there is a
total of two critical hyperplanes sitting inside of $\Delta_{i}(p_{i})$: one is the vertical hyperplane and the other is the
tangency hyperplane, labeled by the letter $T$.  Now, figure $3$ describes the picture for $\Delta_{i}(p_{i})$ when the i-th letter
in the $RVT$ code of $p_{i}$ is the letter $L$.  Figure $(c)$ depicts this situation where there is now a total of three 
critical hyperplanes: one for the vertical hyperplane, and two tangency hyperplanes, labeled as $T_{1}$ and $T_{2}$.  Now, because of
the presence of these three critical hyperplanes we need to refine our notion of an $L$ direction and add two more distinct $L$ directions.
These three directions are labeled as $L_{1}$, $L_{2}$, and $L_{3}$.
\

With the above in mind, we can now refine our $RC$ coding and define the $RVT$ code for points within the Monster tower.
Take $p_{k} \in \mathcal{P}^{k}$ and if $\omega_{i}(p_{k}) = C$ then we look at the point $p_{i} = \pi_{k,i}(p_{k})$, where
$p_{i} = (p_{i-1}, \ell_{i-1})$.  Then depending on which hyperplane $\ell_{i-1}$ is contained in we relabel the letter $C$ by
the letter $V$, $T$, $L$, $T_{i}$ for $i = 1,2$, or $L_{j}$ for $j = 1,2,3$.  As a result, we see that each of the first four
levels of the Monster tower is made up of the following $RVT$ classes:
\

\begin{itemize}
\item{}Level 1:  $R$. 
\item{}Level 2:  $RR, RV$.
\item{}Level 3:  $$RRR, RRV, RVR, RVV, RVT, RVL$$
\item{}Level 4:  
$$RRRR, RRRV$$
$$RRVR, RRVV, RRVT, RRVL$$
$$RVRR, RVRV,  RVVR, RVVV, RVVT, RVVL $$
$$RVTR, RVTV, RVTT , RVTL$$
$$RVLR, RVLV, RVLT_1, RVLT_2, RVLL_1,  RVLL_2, RVLL_3$$
\end{itemize}
\

\begin{rem}
As was pointed out in \cite{castro}, the symmetries, at any level in the Monster tower, preserve the critical hyperplanes.
In other words, if $\Phi^{k}$ is a symmetry at level $k$ in the Monster tower and $\delta^{j}_{i}$ is a critical hyperplane within
$\Delta_{k}$, then $\Phi^{k}_{\ast}(\delta^{j}_{i}) = \delta^{j}_{i}$.  As a result, the $RVT$ classes creates a partition of the various
points within any level of the Monster tower.
\end{rem}

Now, from the above configurations of critical hyperplanes section one might ask the following question: how does one "see" the two
tangency hyperplanes that appear over an "L" point and where do they come from?  This question was an important one to ask when
trying to classify the number of orbits within the fourth level of the Monster Tower and to better
understand the geometry of the tower.  We will provide an example to answer this question, but
before we do so we must discuss some details about a particular coordinate system called Kumpera-Rubin coordinates to help us do various
computations on the Monster tower.
\

\subsection{Kumpera-Rubin coordinates} When doing local computations in the tower (\ref{eqn:tower}) one needs to work with suitable coordinates.
A good choice of coordinates were suggested by Kumpera and Ruiz \cite{kumpera1} in the Goursat case, and later generalized by Kumpera and Rubin \cite{kumpera2} for multi-flags. A detailed description of the inductive construction of Kumpera-Rubin coordinates was given in \cite{castro} and is discussed in the example following this section, as well as in the proof of our level $3$ classification. For the sake of clarity, we will highlight
the coordinates' attributes through an example.

\begin{eg}[Constructing fiber affine coordinates in  $\mathcal{P}^{2}$]
\end{eg}
\subsection*{Level One:}
Consider the pair $(\mathbb{R}^{3}, T \R^{3})$ and let $(x,y,z)$ be local coordinates. The set of 1-forms $\{dx,dy,dz\}$ form a coframe of $T^*\mathbb{R}^3$. 
Any line $\ell$ through $p \in \mathbb{R}^3$ has projective coordinates $[dx(\ell): dy(\ell): dz(\ell)]$. 
Since the affine group, which is contained in $Diff(3)$, acts transitively on $\mathcal{P}(T\mathbb{R}^3)$ we can fix $p=(0,0,0)$ and $\ell = \text{span} \left\{ (1,0,0) \right\}$. 
Thus $dx(\ell)\neq 0$ and we introduce fiber affine coordinates $[1: dy/dx: dz/dx]$ or, $$u = \frac{dy}{dx}, v = \frac{dz}{dx}.$$ 
The Pfaffian system describing the prolonged distribution $\Delta_1$ on $$\mathcal{P}^{1} \approx \mathbb{R}^3\times \mathbb{P}^2$$ is 
$$ \{dy - u dx = 0, dz - v dx = 0 \} = \Delta_1 \subset T\mathcal{P}^{1}.$$

At the point $p_{1} = (p_{0}, \ell) = (x,y,z,u,v) = (0,0,0,0,0)$ the distribution is the linear subspace $$\Delta_1(0,0,0)=\{dy = 0,dz = 0\}.$$ The triad of one-forms $dx,du,dv$ form a local coframe for $\Delta_1$ near $p_{1} = (p_{0}, \ell)$. 
The fiber, $F_{1}(p_{1}) = \pi^{-1}_{1,0}(p_{0})$, is given by $x = y = z = 0$.  The 2-plane of critical directions (``bad-directions'') is thus spanned by $\frac{\partial}{\partial u},\frac{\partial}{\partial v}$.

\

The reader may have noticed that we could have instead chosen any regular direction at level $1$, instead, e.g. $\frac{\partial}{\partial x} + a \frac{\partial}{\partial u} + b \frac{\partial}{\partial v}$ and centered our chart on it.  All regular directions, at level one, are equivalent.

\subsection*{Level Two: $RV$ points.}
Any line $\ell\subset \Delta_1(p_{1}')$, for $p_{1}'$ near $p_{1}$, will have projective coordinates 
$$[dx(\ell): du(\ell): dv(\ell)] .$$
If we choose a critical direction, say $\ell =\frac{\partial}{\partial u}$, then $du(\frac{\partial}{\partial u})=1$ and we can take the projective chart $[\frac{dx}{du} : 1 : \frac{dv}{du}]$.  We will show below that any two critical directions are equivalent and therefore such a choice does not  result in any loss of generality. We introduce new fiber affine coordinates 
$$u_2 = \frac{dx}{du},v_2 = \frac{dy}{du},$$
and the distribution $\Delta_2$ will be described in this chart as 

\begin{eqnarray*}
\Delta_2 = \{dy - u dx = 0, dz - v dx = 0,\\
dx - u_2 du = 0, dv - v_2 du = 0\} \subset T\mathcal{P}^{2}.
\end{eqnarray*}

\subsection*{Level Three: The Tangency Hyperplanes over an $L$ point.}
We take $p_{3} = (p_{2}, \ell) \in RVL$ with $p_{2}$ as in the level two discussion.
We will show what the local affine coordinates near this point are and that the tangency hyperplane $T_{1}$, in $\Delta_{3}(p_{3})$, is the critical hyperplane $\delta^{1}_{2}(p_{3}) = 
span \{ \frac{\pa}{\pa v_{2}}, \frac{\pa}{\pa v_{3}} \}$ and the tangency hyperplane $T_{2}$ is the critical hyperplane
$\delta^{2}_{1}(p_{3}) = span \{ \frac{\pa}{\pa v_{2}}, \frac{\pa}{\pa u_{3}}\}$.
\

We begin with the local coordinates near $p_{3}$.  First, recall that the distribution $\Delta_{2}$, in this case, is coframed
by $[du: du_{2}: dv_{2}]$.  Within $\Delta_{2}$ the vertical hyperplane is given by $du = 0$ and the tangency hyperplane by 
$du_{2} = 0$.  The point $p_{3} = (p_{2}, \ell)$ with $\ell$ being an $L$ direction means that both $du(\ell) = 0$
and $du_{2}(\ell) = 0$.  This means that the only choice for local coordinates near $p_{3}$ is given by
$[\frac{du}{dv_{2}}: \frac{du_{2}}{dv_{2}}: 1]$ to give that the fiber coordinates at level $3$ are 
$$u_{3} = \frac{du}{dv_{2}}, v_{3} = \frac{du_{2}}{dv_{2}} $$ 
and the distribution $\Delta_{3}$ will be described in this chart as 

\begin{eqnarray*}
\Delta_{3} = \{dy - u dx = 0, dz - v dx = 0,\\
dx - u_2 du = 0, dv - v_2 du = 0, \\
du - u_{3}dv_{2} = 0, du_{2} - v_{3}dv_{2} = 0 \} \subset T\mathcal{P}^{3}.
\end{eqnarray*} 

With this in mind, we are ready to determine how the two tangency hyperplanes are situated within $\Delta_{3}$.

\ 

{\it $T_{1} = \delta^{1}_{2}(p_{3})$:}
First we note that $p_{3} = (x,y,z,u,v,u_{2},v_{2},u_{3},v_{3}) = (0,0,0,0,0,0,0,0,0)$ with $u = \frac{dy}{dx}$, 
$v = \frac{dz}{dx}$, $u_{2} = \frac{dx}{du}$, $v_{2} = \frac{dv}{du}$, $u_{3} = \frac{du}{dv_{2}}$, 
$v_{3} = \frac{du_{2}}{dv_{2}}$.  With this in mind, we start by looking at the vertical hyperplane 
$V_{2}(p_{2}) \subset \Delta_{2}(p_{2})$ and prolong the fiber $F_{2}(p_{2})$ associated to $V_{2}(p_{2})$ and see that

$$ \mathcal{P}^{1}(F_{2}(p_{2})) = \Pj V_{2} = (p_{1}, u_{2}, v_{2}, [ du: du_{2}: dv_{2} ]) = 
(p_{1}, u_{2}, v_{2}, [ 0: a: b ]) $$ 
$$= (p_{1}, u_{2}, v_{2}, [ 0: \frac{a}{b}: 1 ]) = (p_{1}, u_{2}, v_{2}, 0, v_{3})$$
where $a,b \in \R$ with $b \neq 0$.  Then, since $\Delta_{3}$, in a neighborhood of $p_{3}$, is given by 
$$ \Delta_{3} = span \{ u_{3}Z^{(2)}_{1} + v_{3}\frac{\pa}{\pa u_{2}} + \frac{\pa}{\pa v_{2}}, \frac{\pa}{\pa u_{3}}, \frac{\pa}{\pa v_{3}} \}$$
with $Z^{(2)}_{1} = u_{2}Z^{(1)}_{1} + \frac{\pa}{\pa u} + v_{2} \frac{\pa}{\pa v}$ and 
$Z^{(1)}_{1} = u \frac{\pa}{\pa y} + v \frac{\pa}{\pa z} + \frac{\pa}{\pa x}$ and that 
$T_{p_{3}}(\mathcal{P}^{1}(F_{2}(p_{2}))) = span \{ \frac{\pa}{\pa u_{2}}, \frac{\pa}{\pa v_{2}}, \frac{\pa}{\pa v_{3}} \}$ we see that since
$$\delta^{1}_{2}(p_{3}) = \Delta_{3}(p_{3}) \cap T_{p_{3}}(\mathcal{P}^{1}(F_{2}(p_{2})))$$
it gives us that
$$\delta^{1}_{2}(p_{3}) = span \{ \frac{\pa}{\pa v_{2}}, \frac{\pa}{\pa v_{3}} \}$$
Now, since $V_{3}(p_{3}) \subset \Delta_{3}(p_{3})$ is given by $V_{3}(p_{3}) = span \{ \frac{\pa}{\pa u_{3}}, \frac{\pa}{\pa v_{3}} \}$    
we see, based upon figure $(c)$, that $T_{1} = \delta^{1}_{2}(p_{3})$.

{\it $T_{2} = \delta^{2}_{1}(p_{3})$:}
We begin by looking at $V_{1}(p_{1}) \subset \Delta_{1}(p_{1})$ and at the fiber $F_{1}(p_{1})$ associated to $V_{1}(p_{1})$.
When we prolong the fiber space we see that  
$$\mathcal{P}^{1}(F_{1}(p_{1})) = \Pj V_{1} = (0,0,0, u, v, [dx: du: dv]) = (0,0,0, u, v, [0: a: b])$$
$$ = (0,0,0, u, v, [0: 1: \frac{b}{a}]) = (0,0,0, u, v, 0, v_{2})$$
where $a,b \in \R$ with $a \neq 0$.  Then, since $\Delta_{2}$, in a neighborhood of $p_{2}$, is given by
$$ \Delta_{2} = span \{ u_{2}Z^{(1)}_{1} + \frac{\pa}{\pa u} + v_{2}\frac{\pa}{\pa v}, \frac{\pa}{\pa u_{2}}, \frac{\pa}{\pa v_{2}} \}$$
and $T_{p_{2}}(\mathcal{P}^{1}(F_{1}(p_{1}))) = span \{ \frac{\pa}{\pa u}, \frac{\pa}{\pa v}, \frac{\pa}{\pa v_{2}} \}$ that
$$\delta^{1}_{1}(p_{2}) = \Delta_{2}(p_{2}) \cap T_{p_{2}}(\mathcal{P}^{1}(F_{1}(p_{1})))$$
and we see that in a neighborhood of $p_{2}$ that 
$$\delta^{1}_{1} = span \{ u_{2}Z^{(1)}_{1} + \frac{\pa}{\pa u} + v_{2}\frac{\pa}{\pa v}, \frac{\pa}{\pa v_{2}} \}$$
Now, in order to figure out what $\delta^{2}_{1}(p_{3})$ is we need to prolong the fiber $F_{1}(p_{1})$ twice and then look at
the tangent space at the point $p_{3}$.  We see that 
$$\mathcal{P}^{2}(F_{1}(p_{1})) = \Pj \delta^{1}_{1} =  (0,0,0, u,v, 0, v_{2}, [du: du_{2}: dv_{2}])$$
$$ = (0,0,0, u,v, 0, v_{2}, [a: 0: b]) = (0,0,0,u, v, 0, v_{2}, [\frac{a}{b}: 0: 1]) = (0,0,0, u, v, 0, v_{2}, u_{3}, 0)$$
then since 
$$\delta^{2}_{1}(p_{3}) = \Delta_{3}(p_{3}) \cap T_{p_{3}}(\mathcal{P}^{2}(F_{1}(p_{1})))$$
with $\Delta_{3}(p_{3}) = span \{ \frac{\pa}{\pa v_{2}}, \frac{\pa}{\pa u_{3}}, \frac{\pa}{\pa v_{3}} \}$ and
$T_{p_{3}}(\mathcal{P}^{2}(F_{1}(p_{1}))) = span \{ \frac{\pa}{\pa u}, \frac{\pa}{\pa v}, \frac{\pa}{\pa v_{2}}, \frac{\pa}{\pa u_{3}} \}$
it gives that
$$\delta^{2}_{1}(p_{3}) = span \{ \frac{\pa}{\pa v_{2}}, \frac{\pa}{\pa u_{3}} \}$$
and from looking at figure $(c)$ one can see that $T_{2} = \delta^{2}_{1}(p_{3})$.

%$F_{3}$, $\mathcal{P}^{3}$, $\mathcal{P}^{2}$, $\R^{3}$, $\pi$, $F_{2}$, $F_{1}$, 
%$\mathcal{P}^{1}(F_{2}(p_{2}))$, $\mathcal{P}^{1}(F_{1}(p_{1}))$, $p_{1}$, $p_{2}$, $p_{3}$
%$\frac{\pa}{\pa u_{3}}$,  $\frac{\pa}{\pa u_{3}} + \frac{\pa}{\pa v_{3}}$, $\Delta_{3}$
%$\R^{3}$, $u, v$, $\Delta_{0}$, $\ell$, $\Delta_{1} = d \pi^{-1}(\ell)$, $p$, $q$, $\Phi$, 
%$\Phi^{1}$, $p_{1}$, $q_{1}$

\begin{figure}
  \centering
  \includegraphics[width=0.4\textwidth]{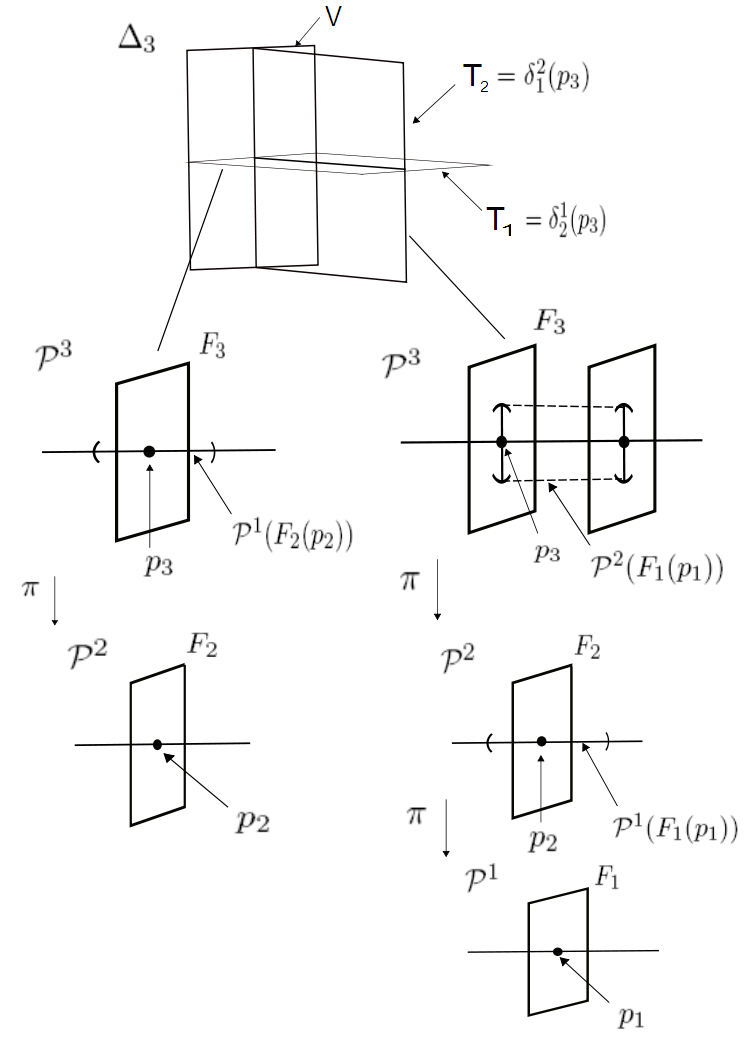}              
  \caption{Critical hyperplane configuration over $p_{3} \in RVL$.}
\end{figure}

\begin{rem}
The above example, along with figure $2$, gives concrete reasoning for why a critical hyperplane, which is not the vertical one,
is called a "tangency" hyperplane.  Also, in figure $2$ we have drawn the submanifolds $\mathcal{P}^{1}(F_{2}(p_{2}))$ and
$\mathcal{P}^{1}(F_{1}(p_{1}))$ to reflect the fact that they are tangent to the manifolds $\mathcal{P}^{3}$ and $\mathcal{P}^{2}$
respectively with one of their dimensions tangent to the vertical space.  At the same time, the submanifold 
$\mathcal{P}^{2}(F_{1}(p_{1}))$ is drawn to reflect that it is tangent to the manifold $\mathcal{P}^{3}$ with one of its dimensions
tangent to the appropriate direction in the vertical hyperplane.  In particular, it is drawn to show the fact that 
$\mathcal{P}^{2}(F_{1}(p_{1}))$ is tangent to the $\frac{\partial}{\partial u_{3}}$ direction while 
$\mathcal{P}^{1}(F_{2}(p_{2}))$ is tangent to the $\frac{\partial}{\partial v_{3}}$ direction.  
\end{rem} 

\subsection{Semigroup of a curve}

\begin{defn}
The \textit{order} of an analytic curve germ $f(t) = \Sigma a_{i}t^{i}$ is the smallest integer $i$ such that 
$a_{i} \neq 0$.  We write $ord(f)$ for this (nonegative) integer.  The \textit{multiplicity} of a curve germ
$\gamma:(\mathbb{R},0) \to (\mathbb{R}^n,0)$, denoted $mult(\gamma)$, is the minimum of the orders of its coordinate functions 
$\gamma_{i}(t)$ relative to any coordinate system vanishing at $p$. 
\end{defn}

\begin{defn}
If $\gamma: (\mathbb{R},0) \to (\mathbb{R}^n,0)$ is a well-parameterized curve germ, then
its \emph{semigroup} is the collection of positive integers $\text{ord}(P(\gamma(t)))$ as
$P$ varies over analytic functions of $n$ variables vanishing at $0$.
\end{defn}  Because $\text{ord}(PQ(\gamma(t))) = \text{ord}(P(\gamma(t)) + \text{ord}(Q(\gamma(t))$
the curve semigroup is indeed an algebraic semigroup, i.e. a subset of $\mathbb{N}$ closed under addition.
The semigroup of a well-parameterized curve is a basic  diffeomorphism invariant of the curve.    
\begin{defn}[Following Arnol'd (cf.\cite{arnsing}, end of his introduction] 
A curve germ $\gamma$ in $\mathbb{R}^d$, $d \ge 2$  
has symbol  $[m, n]$,  $[m, n, p]$, or $[m, (n, p)]$  if it is equivalent to a  curve germ of the form 
$(t^m, t^n,  0,  \ldots,0 ) + O(t^{n+1})$,  
$(t^m, t^n, t^p, 0, \ldots,0) + O(t^{p+1})$, or
$(t^m, t^n + t^p, 0, \ldots,0) + O(t^{p+1})$. 
Here $m< n < p$ are positive integers.  
\end{defn}

\subsection{The points-to-curves and back philosophy}

The idea is to translate the problem of classifying orbits in the tower (\ref{eqn:tower}) into an equivalent classification problem for finite jets of space curves. Here we are going to mention some highlights of this approach, we will refer the diligent reader to \cite{castro} check the technical details.
\subsubsection*{Methodology. How does the curve-to-point philosophy work?} 
To any $p\in \mathcal{P}^{k}(n)$ we associate the set 
\begin{eqnarray*}
\text{Germ}(p) &=& \{ c :(\mathbb{R},0)\rightarrow (\mathbb{R}^3,0)| \text{$\frac{dc^{k}}{dt}\vert_{t=0}\neq 0$ is a regular direction} \}.
\end{eqnarray*}
The operation of $k$-fold prolongation applied to $\text{Germ}(p)$ yields immersed curves at level $k$ in the Monster tower, and tangent to some line $\ell$ having nontrivial projection onto the base manifold $\mathbb{R}^3.$ Such ``good directions'' were named {\it regular} in \cite{castro} and within each subspace $\Delta_k$ they form an open dense set.  
A ``bad direction'' $\ell_{\text{critical}}$, or {\it critical direction} in the jargon of \cite{castro} are directions which will project down to a point. The set of critical direction has codimension 1, and consists of a finite union of 2-planes within each $\Delta_k$. Symmetries of $\mathcal{P}^{k}$ do preserve the different types of directions. 

In \cite{castro} it was proved that that $\text{Germ}(p)$ is always non-empty. 
Consider now the set valued map $p\mapsto \text{Germ}(p)$. 
One can prove that $p\sim q$ iff $\text{Germ}(p)\sim \text{Germ}(q)$. (The latter equivalence means ``to any curve in $\text{Germ}(p)$ there is a curve in $\text{Germ}(q)$ and vice-versa.'') 
\begin{lemma}[Fundamental lemma of points-to-curves approach] 
Let $\Omega$ be a subset of $\mathcal{P}^{k}(n)$ and suppose that 
$$\bigcup_{p\in \Omega}\text{Germ}(p)= \{\text{finite no. of equivalence classes of curves}\}.$$
Then $$\Omega = \{\text{finite no. of orbits}\}.$$ 
\end{lemma}

\section{Proofs.}
Now we are ready to prove Theorem \ref{thm:main} and Theorem \ref{thm:count}.  We start at level $1$ of the tower and work our way
up to level $4$.  At each level we classify the orbits within each of the $RVT$ classes that can arise at that
particular level.  We begin with the first level.

\noindent {\bf Proof of Theorem \ref{thm:main} and Theorem \ref{thm:count}, the classification of points at level $1$ and level $2$.}
Theorem $3.3$ tells us that all points at the first level of the tower are equivalent, giving that there is a single 
orbit.  For level $2$ there are only two possible $RVT$ codes: $RR$ and $RV$.  Again, any point in the class $RR$ is 
a Cartan point and by Theorem $3.4$ consists of only one orbit.  The class $RV$ consists of a single orbit by
Theorem $3.5$.

\

{\bf The classifiction of points at level $3$.}
There is a total of six distinct $RVT$ classes at level three in the Monster tower.  We begin with the class $RRR$.

{\it The class $RRR$:}
Any point within the class $RRR$ is a Cartan point that by Theorem \ref{thm:cartan} that there is only one orbit within
this class.

{\it The classes $RVR$ and $RRV$:}
From Theorem \ref{thm:ak} we know that any point within the class $RVR$ has a single orbit, which
is represented by the point $\gamma^{3}(0)$ where $\gamma$ is the curve $\gamma(t) = (t^{2}, t^{3}, 0)$.
Similarly, the class $RRV$ has a single orbit, which is represented by the point $\tilde{\gamma}^{3}(0)$ where
$\tilde{\gamma}(t) = (t^{2}, t^{5}, 0)$. 

\no Before we continue, we need to pause and provide some framework to help us with the classification of the
remaining $RVT$ codes.

\

{\it Setup for classes of the form $RVC$:}
We set up coordinates $x,y,z,u,v,u_{2},v_{2}$ for a point in the class $RV$ as in subsection $4.4$.
Then for $p_{2} \in RV$ we have $\Delta_{2}(p_{2}) = span \{ \frac{\pa}{\pa u}, \frac{\pa}{\pa u_{2}}, \frac{\pa}{\pa v_{2}} \}$
where $p_{2} = (x,y,z,u,v,u_{2}, v_{2}) = (0,0,0,0,0,0,0)$ and for any point $p_{3} \in RVC \subset \mathcal{P}^{3}$ that
$p_{3} = (p_{2}, \ell_{2}) = (p_{2}, [du(\ell_{2}): du_{2}(\ell_{2}): dv_{2}(\ell_{2})])$.  Since the point $p_{2}$ is in the class $RV$ we see that
if $du = 0$ along $\ell_{2}$ that $p_{3} \in RVV$, $du_{2} = 0$ with $du \neq 0$ along $\ell_{2}$ that $p_{3} \in RVT$,
and if $du = 0$ and $du_{2} = 0$ along $\ell_{2}$ that $p_{3} \in RVL$.  With this in mind, we are ready to continue 
with the classification. 

\

{\it The class $RVV$:}
Let $p_{3} \in RVV$ and let $\gamma \in Germ(p_{3})$.  We prolong $\gamma$ three times and have that 
$\gamma^{3}(t) = (x(t),y(t),z(t), u(t), v(t), u_{2}(t), v_{2}(t))$ and we look at the component functions 
$u(t)$, $u_{2}(t)$, and $v_{2}(t)$ where we set $u(t) = \Sigma_{i} a_{i}t^{i}$, $u_{2}(t) = \Sigma_{j}b_{j}t^{j}$, and 
$v_{2}(t) = \Sigma_{k}c_{k}t^{k}$.  Now, since $\gamma^{2}(t)$ needs to be tangent to the vertical hyperplane in $\Delta_{3}$
that $\gamma^{2}(0)'$ must be a proper vertical direction in $\Delta_{3}$, meaning $\gamma^{2}(0)'$ is not an $L$ direction.
Since $\Delta_{3}$ is coframed by $du$, $du_{2}$, and $dv_{2}$, we must have that $du = 0$ and $du_{2} \neq 0$ along
$\gamma^{2}(0)'$.  This imposes the condition for the functions $u(t)$ and $u_{2}(t)$ that $a_{1} = 0$ and $b_{1} \neq 0$, but
for $v_{2}(t)$ it may or may not be true that $c_{1}$ is nonzero.  Also it must be true that $a_{2} \neq 0$ or else the curve
$\gamma$ will not be in the set $Germ(p_{3})$.  We first look at the case when $c_{1} \neq 0$. 

{\it Case $1$, $c_{1} \neq 0$:}
From looking at the one-forms that determine $\Delta_{2}$, we see that in order for the curve $\gamma^{3}$ to be integral
to this distribution that the other component function for $\gamma^{3}$ must satisfy the following relations:
$$ \dot{y}(t) = u(t) \dot{x}(t), \dot{z}(t) = v(t) \dot{x}(t)$$
$$ \dot{x}(t) = u_{2}(t) \dot{u}(t), \dot{v}(t) = v_{2}(t) \dot{u}(t)$$
We start with the expressions for $\dot{x}(t)$ and $\dot{v}(t)$ and see that, based upon what we know about
$u(t)$, $u_{2}(t)$, and $v_{2}(t)$, that $x(t) = \frac{2a_{2}b_{1}}{3}t^{3} + \ldots$ and 
$v(t) = \frac{2a_{2}c_{1}}{3}t^{3} + \ldots$.  We can then use this information to help us find $y(t)$ and $z(t)$.
We see that $y(t) = \frac{2a^{2}_{2}b_{1}}{5}t^{5} + \ldots$ and $z(t) = \frac{4a^{2}_{2}b_{1}c_{1}}{3}t^{7} + \ldots$. 
Now, we know what the first nonvanishing coeffiecients are for the curve $\gamma(t) = (x(t), y(t), z(t))$ and we want to
determine the simpliest curve that $\gamma$ must be equivalent to.  In order to do this we will first look at the semigroup
for the curve $\gamma$.  In this case the semigroup is given by $S = \{3, [4], 5, 6, 7, \cdots \}$.

\begin{rem}
We again pause to explain the notation used for the semigroup $S$.  The set $S = \{3, [4], 5, 6, 7, \cdots \}$ is a semigroup
where the binary operation is addition.  The numbers $3$, $5$, $6$, and so on are elements of this semigroup while the bracket
around the number $4$ means that it is not an element of $S$.  When we write "$\cdots$" after the number $7$ it means that 
every positive integer after $7$ is an element in our semigroup.
\end{rem}

This means that every term, $t^{i}$ for $i \geq 7$, can be eliminated
from the above power series expansion for the component functions $x(t)$, $y(t)$, and $z(t)$ by a change of variables
given by $(x,y,z) \mapsto (x + f(x,y,z), y + g(x,y,z), z + h(x,y,z))$.  With this in mind and 
after we rescale the leading coeffiecients for each of the component for $\gamma$, we see that 
$$\gamma(t) = (x(t), y(t), z(t)) \sim (\tilde{x}(t), \tilde{y}(t), \tilde{z}(t)) = (t^{3} + \alpha t^{4}, t^{5}, t^{7}) $$
We now want to see if we can eliminate the $\alpha$ term, if it is nonzero.  To do this we will use a combination of reparametrization
techniques along with semigroup arguments.  Use the reparamentrization $t = T(1 - \frac{\alpha}{3}T)$ and we get that
$\tilde{x}(T) = T^{3}(1 - \frac{\alpha}{3}T)^{3} + T^{4}(1 - \frac{\alpha}{3}T)^{4} + \ldots = T^{3} + O(T^{5})$.
This gives us that $(\tilde{x}(T), \tilde{y}(T), \tilde{z}(T)) = (T^{3} + O(T^{5}), T^{5} + O(T^{6}), T^{7} + O(T^{8}))$
and since we can eliminate all of the terms of degree $5$ and higher we see that
$(\tilde{x}(T), \tilde{y}(T), \tilde{z}(T)) \sim (T^{3}, T^{5}, T^{7})$.  This means that our original $\gamma$ is equivalent
to the curve $(t^{3}, t^{5}, t^{7})$.

{\it Case $2$, $c_{1} = 0$:}
By repeating an argument similar to the above one, we will end up with $\gamma(t) = (x(t), y(t), z(t)) = 
(\frac{2a_{2}b_{1}}{3}t^{3} + \ldots, \frac{2a^{2}_{2}b_{1}}{5}t^{5} + \ldots, \frac{a^{2}_{2}b_{1}c_{2}}{8}t^{8} + \ldots)$.
Note that $c_{2}$ may or may not be equal to zero.  This gives that the semigroup for the curve $\gamma$ is $S = \{ 3, [4], 5, 6, [7], 8  \cdots \}$ and that our curve $\gamma$ is such that
$$\gamma(t) = (x(t), y(t), z(t)) \sim (\tilde{x}(t), \tilde{y}(t), \tilde{z}(t)) = 
(t^{3} + \alpha_{1}t^{4} + \alpha_{2}t^{7} ,t^{5} + \beta t^{7}, 0)$$
Again, we want to know if we can eliminate the $\alpha_{i}$ and $\beta$ terms.  First we focus on the $\alpha_{i}$ terms
in $\tilde{x}(t)$.  We use the reparametrization given by $t = T(1 - \frac{\alpha_{1}}{3}T)$ to give us that
$\tilde{x}(T) = T^{3} + \alpha_{2}'T^{7} + O(T^{8})$.  Then to eliminate the $\alpha_{2}'$ term we use the reparametrization
given by $T = S(1 - \frac{\alpha_{2}'}{3}S^{4})$ to give that $\tilde{x}(S) = S^{3} + O(S^{8})$.
We are now ready to deal with the $\tilde{y}$ function.  Now, because of our two reparametrizations we get that
$\tilde{y}$ is of the form $\tilde{y}(t) = t^{5} + \beta't^{7}$.  To get rid of the $\beta'$ term we simply use the 
rescaling given by $t \mapsto \frac{1}{ \sqrt{ \left|\beta'  \right|} }t$ and then use the scaling diffeomorphism
given by $(x,y,z) \mapsto ( \left| \beta' \right|^{\frac{3}{2}}x, \left| \beta' \right|^{\frac{5}{2}}y, z)$
to give us that $\gamma$ is equivalent to either $(t^{3}, t^{5} + t^{7}, 0)$ or $(t^{3},t^{5} - t^{7}, 0)$.
Note that the above
calculations were done under the assumption that $\beta_{1} \neq 0$.  If $\beta_{1} = 0$ then we see, using similar calculations
as the above, that we get the normal form $(t^{3}, t^{5}, 0)$.  This means that there is a total of $4$ possible normal forms
that represent the points within the class $RVV$.  It is tempting, at first glance, to believe that these curves are all inequivalent.
However, it can be shown that the $3$ curves $(t^{3}, t^{5} + t^{7}, 0)$, $(t^{3},t^{5} - t^{7}, 0)$, $(t^{3}, t^{5}, 0)$
are actually equivalent.  It is not very difficult to show this equivalence, but it does amount to rather messy calculation.
As a result, the techniques used to show this equivalence are outlined in section $7.1$ of the appendix.
\

This means that the total number of possible normal forms is reduced to $2$ possibilities: $\gamma_{1}(t) = (t^{3}, t^{5}, t^{7})$ and 
$\gamma_{2}(t) = (t^{3}, t^{5}, 0)$.
We will show that these two curves are inequivalent to one another.  One possibility is to look at the semigroups that each of these
curves generate.  The curve $\gamma_{1}$ has the semigroup $S_{1} = \{3,[4], 5, 6, 7, \cdots \}$, while the curve $\gamma_{2}$ has the semigroup
$S_{2} = \{3, [4], 5, 6, [7], 8, \cdots \}$.  Since the semigroup of a curve is an invariant of the curve and the two curves generate different
semigroups, then the two curves must be inequivalent.  In \cite{castro} another method was outlined to check and see whether
or not these two curves are equivalent, which we will now present.  One can see that the curve $(t^{3}, t^{5}, 0)$ is a planar curve and 
in order for the curve $\gamma_{1}$ to be equivalent to the curve $\gamma_{2}$, then we must be able to find a way to turn 
$\gamma_{1}$ into a planar curve, meaning we need to find a change of variables or a reparametrization which 
will make the $z$-component function of $\gamma_{1}$ zero.  If it were true that $\gamma_{1}$ is actually a planar curve, then $\gamma_{1}$ 
must lie in an embedded surface in $\R^{3}$(or embedded surface germ), say $M$.  Since $M$ is an embedded surface it means that there exists a local defining
function at each point on the manifold.  Let the local defining function near the origin be the real analytic function $f: \R^{3} \rightarrow \R$.  Since $\gamma_{1}$ is on M, then $f(\gamma_{1}(t)) = 0$ for all $t$ near zero.
However, when one looks at the individual terms in the Taylor series expansion of $f$ composed with $\gamma_{1}$ that there will
be nonzero terms which will show up and give that $f(\gamma_{1}(t)) \neq 0$ for all $t$ near zero, which creates a contradiction.  This tell us that $\gamma_{1}$ cannot be equivalent 
to any planar curve near $t=0$.  As a result, there is a total of two inequivalent normal forms for 
the class $RVV$: $(t^{3}, t^{5}, t^{7})$ and $(t^{3}, t^{5}, 0)$.
\

The remaining classes $RVT$ and $RVL$ are proved in an almost identical manner using the above ideas and techniques.  As a result,
we will omit the proofs and leave them to the reader.  
\

With this in mind, we are now ready to move on to the fourth level of the tower.  We initially tried to tackle the problem of classifying the orbits
at the fourth level by using the curve approach from the third level.  Unfortunately, the curve approach becomes a bit
too unwieldy to use to determine what the normal forms were for the various $RVT$ classes.  
The problem was simply this: when we looked at the semigroup
for a particular curve in a number of the $RVT$ classes at the fourth level that there were too many "gaps" in the various
semigroups. The first occuring class, according to codimension, in which this occured was the class $RVVV$. 

\begin{eg}  
The semigroups for the class $RVVV$.  Let $p_{4} \in RVVV$, and for $\gamma \in Germ(p_{4})$ that 
$\gamma^{3}(t) = (x(t), y(t), z(t), u(t), v(t), u_{2}(t), v_{2}(t), u_{3}(t), v_{3}(t))$ with 
$u = \frac{dy}{dx}$, $v = \frac{dz}{dx}$, $u_{2} = \frac{dx}{du}$, $v_{2} = \frac{dv}{du}$, $u_{3} = \frac{du}{du_{2}}$, 
$v_{3} = \frac{dv_{2}}{du_{2}}$. Since $\gamma^{4}(0) = p_{4}$ we must have that $\gamma^{3}(t)$ is tangent to the vertical
hyperplane within $\Delta_{3}$, which is coframed by $\left\{ du_{2}, du_{3}, dv_{3} \right\}$.  One can see that $du_{2} = 0$ along 
$\gamma^{3}(0)'$.  Then, just as with the analysis with the third level, we look at $u_{2}(t) = \Sigma_{i} a_{i}t^{i}$,
$u_{3}(t) = \Sigma_{j}b_{j}t^{j}$, $v_{3}(t) = \Sigma_{k}c_{k}t^{k}$ where we must have that $a_{1} = 0$, $a_{2} \neq 0$,
$b_{1} \neq 0$, and $c_{1}$ may or may not be equal to zero.  When we go from the fourth level back down to the zeroth level
we see that $\gamma(t) = (t^{5} + O(t^{11}), t^{8} + O(t^{11}), O(t^{11}))$.  If $c_{1} \neq 0$, then we get that
$\gamma(t) = (t^{5} + O(t^{12}), t^{8} + O(t^{12}), t^{11} + O(t^{12}))$ and the semigroup for this curve is
$S = \{5, [6], [7], 8, [9], 10, 11, [12], 13, [14], 15, 16, [17], 18 \cdots \}$.  If $c_{1} \neq 0$, then we get that
$\gamma(t) = (t^{5} + O(t^{12}), t^{8} + O(t^{12}), O(t^{12}))$ and the semigroup for this curve is
$S = \{5, [6], [7], 8, [9], 10, [11], [12], 13, [14], 15, 16, [17], 18, [19], 20, 21, [22], 23 \cdots \}$.   

\end{eg}

As a result, it became impractical to work stictly using the curve approach.
This meant that we had to look at a different approach to the classification problem.  This lead us to work with a tool
called the isotropy representation.

\subsection{The isotropy method.}

Here is the general idea of the method. 

Suppose we want to look at a particular RVT class, at the k-th level, given by $\omega_k$ (a word of
length $k$) and we want to see how many orbits there are.  Suppose as well that we understand its projection 
$\omega_{k-1}$ one level down, which decomposes into $N$ orbits.  Choose representative points $p_{i}$, $i = 1, \cdots , N$ for
the $N$ orbits in $\omega_{k-1}$, and consider the group $G_{k-1}(p_{i})$ of level $k-1$ symmetries that fix $p_{i}$.
This group is called the \textit{isotropy group of} $p_{i}$.  Since elements $\Phi^{k-1}$ of the isotropy group
fix $p_{i}$, their prolongations $\Phi^{k} = (\Phi^{k-1}, \Phi^{k-1}_{\ast})$ act on the fiber over $p_{i}$.  
Under the action of the isotropy group the fiber decomposes into some
number $n_{i} \geq 1$ (possibly infinite) of orbits.  Summing, we find that $\omega_{k}$ decomposes into 
$\sum_{i = 1}^{N} n_{i} \geq N$ orbits.
For the record, $\Phi_*^{\bullet}$ denotes the tangent map. 
This will tell us how many
orbits there are for the class $\omega_k$. 

This is the theory. Now we need to explain how one actually prolongs diffeomorphisms in practice. Since the manifold $\mathcal{P}^{k}$ is a type of fiber compactification of $J^k(\mathbb{R},\mathbb{R}^2)$, it is reasonable to expect that the prolongation of diffeomorphisms from the base $\mathbb{R}^3$ should be similar to what one does when prolonging point symmetries from the theory of jet spaces. See specially (\cite{duzhin}, last chapter) and (\cite{olver}, p. 100). 

Given a point $p_k \in \mathcal{P}^{k}$ and a map $\Phi\in \text{Diff}(3)$ we would like to write explicit formulas for 
$$\Phi^{k}(p_k).$$ 
Coordinates of $p_k$ can be made explicit. 
Now take any curve $\gamma(t) \in \text{Germ}(p_k)$, and consider the prolongation of $\Phi\circ \gamma(t)$. The coordinates of $\Phi^{k}(p_k)$ are exactly the coordinates of $(\Phi\circ \gamma)^{(k)}(0) = \Phi^{k}(\gamma^{k}(0))$. Moreover the resulting point is independent of the choice of a regular $\gamma \in \text{Germ}(p)$.

\subsection{Proof of theorem \ref{thm:main}}
\

{\bf The classifiction of points at level $4$.}
We are now ready to begin with the classification of points at level $4$.  We will present the proof for the classification of 
the class $RVVV$ as an example of how the isotropy representation method works.

{\it The class $RVVV$.}
Before we get started, we will summarize the main idea of the following calculation to classify the number of orbits within the
class $RVVV$.  Let $p_{4} \in RVVV \subset \mathcal{P}^{4}$ and start with the projection of $p_{4}$ to level zero,
$\pi_{4,0}(p_{4}) = p_{0}$.  Since all of the points at level zero are equivalent, then one is free to choose any representative
for $p_{0}$. For simplicity, it is easiest to choose it to be the point $p_{0} = (0,0,0)$.  Next, we look at all of the points
at the first level, which project to $p_{0}$.  Since all of these points are equivalent it gives that there is a single orbit in the
first level and we are again able to choose any point in $\mathcal{P}^{1}$ as our representive so long as it projects to the point $p_{0}$.
We will pick $p_{1} = (0,0,0,[1: 0: 0]) = (0,0,0,0,0)$ with $u = \frac{dy}{dx}$ and $v = \frac{dz}{dx}$ and we will look at all of the diffeomorphisms that
fix the point $p_{0}$ and $\Phi_{\ast}([1: 0: 0]) = [1: 0: 0]$.  This condition will place some restrictions on the component functions 
of the local diffeomorphisms $\Phi$ in $Diff_{0}(3)$ when we evaluate at the the point $p_{0}$ and tell us what $\Phi^{1} = (\Phi, \Phi_{\ast})$ will
look like at the point $p_{1}$.  We call this group of diffeomorphisms $G_{1}$.  We can then move on to the second level and look at the class $RV$.  For any $p_{2} \in RV$ it is of the form $p_{2} = (p_{1}, \ell_{1})$ with $\ell_{1}$ contained in the vertical hyperplane inside of $\Delta_{1}(p_{1})$.  Now, apply the pushforwards of the $\Phi^{1}$'s in $G_{1}$ to the vertical hyperplane and see if these symmetries will act transitively on
the critical hyperplane.  If they do act transitively then there is a single orbit within the class $RV$.  If not, then there exists more
than one orbit within the class $RV$.  Note that because of Theorem \ref{thm:ak} that we should expect to only see one orbit within this class.
Once this is done, we can just iterate the above process to classify the number of orbits within the class $RVV$ at the third level and
then within the class $RVVV$ at the fourth level.    

{\it Level 0:} Let $G_{0}$ be the group that contains all diffeomorphism germs that fix the origin.  

{\it Level 1:} We know that all the points in $\mathcal{P}^{1}$ are equivalent, giving that there is only a single orbit.  
So we pick a representative element from the single orbit of $\mathcal{P}^{1}$.  We will take our representative to be $p_{1} = (0,0,0,0,0) = (0,0,0, [1: 0: 0]) = (x,y,z, [dx: dy: dz])$ and 
take $G_{1}$ to be the
set of all $\Phi \in G_{0}$ such that $\Phi^{1}$ will take the tangent to the $x$-Axis back to the 
$x$-Axis, meaning $\Phi_{\ast}([1: 0: 0]) = [1: 0: 0]$.
\

\no Then for $\Phi \in G_{1}$ and $\Phi(x,y,z) = (\phi^{1}, \phi^{2}, \phi^{3})$ we have that 

\

\[
\Phi_{\ast} =
\begin{pmatrix}
  \phi_{x}^{1}  &  \phi_{y}^{1}  &  \phi_{z}^{1} \\
  \phi_{x}^{2}  &  \phi_{y}^{2}  &  \phi_{z}^{2} \\
  \phi_{x}^{3}  &  \phi_{y}^{3}  &  \phi_{z}^{3} \\  
\end{pmatrix}
=
\begin{pmatrix}
  \phi_{x}^{1}  &  \phi_{y}^{1}  &  \phi_{z}^{1} \\
  0             &  \phi_{y}^{2}  &  \phi_{z}^{2} \\
  0             &  \phi_{y}^{3}  &  \phi_{z}^{3} \\
\end{pmatrix}
\]

\

When we evalutate at $(x,y,z) = (0,0,0)$.

\

Here is the \textit{Taylor triangle} representing the different coefficients in the Taylor series expansion of a diffeomorphism
in $G_{i}$.  The three digits represent the number of partial derivatives with respect to either $x$, $y$, or $z$.
For example, $(1,2,0) = \frac{\pa^{3}}{\pa x \pa^{2} y}$.  The vertical column denotes the coefficient order.  We start with the
Taylor triangle for $\phi^{2}$:

\

\begin{center}
\begin{tabular}{rcccccccccc}
$n=0$:&    &    &    &    & \xcancel{(0,0,0)} \\\noalign{\smallskip\smallskip}
$n=1$:&    &    &    &  $\xcancel{(1,0,0)}$ &  (0,1,0) & (0,0,1)\\\noalign{\smallskip\smallskip}
$n=2$:&    &    &  (2,0,0) & (1,1,0)   &  (1,0,1) & (0,1,1)    &  (0,0,2)\\\noalign{\smallskip\smallskip}
\end{tabular}
\end{center}  

\

We have crossed out $(1,0,0)$ since $\frac{\pa \phi^{2}}{\pa y}(\z) = 0$.  Next is the Taylor triangle for
$\phi^{3}$:

\

\begin{center}
\begin{tabular}{rcccccccccc}
$n=0$:&    &    &    &    &  \xcancel{(0,0,0)} \\\noalign{\smallskip\smallskip}
$n=1$:&    &    &    &  $\xcancel{(1,0,0)}$ &  (0,1,0) & (0,0,1)\\\noalign{\smallskip\smallskip}
$n=2$:&    &    &  (2,0,0) & (1,1,0)   &  (1,0,1) & (0,1,1)    &  (0,0,2)\\\noalign{\smallskip\smallskip}
\end{tabular}
\end{center}  

\

This describes the properties of the elements $\Phi \in G_{1}$.

We now try to figure out what the $\Phi^{1}$, for $\Phi \in G_{1}$, will look like in $KR$-coordinates.  First, take $\ell \subset \Delta_{0}$ where we write 
$\ell = a \frac{\partial}{\partial x} + b \frac{\partial}{\partial y} + c \frac{\partial}{\partial z}$ with $a,b,c \in \R$ and
$a \neq 0$.
We see that 

\begin{eqnarray*}
\Phi_{\ast}(\ell)&=&span\left\{ (a \phi_{x}^{1} + b \phi_{y}^{1} + c \phi_{z}^{1})\frac{\pa}{\pa x} + 
 (a \phi_{x}^{2} + b \phi_{y}^{2} + c \phi_{z}^{2})\frac{\pa}{\pa y} +
 (a \phi_{x}^{3} + b \phi_{y}^{3} + c \phi_{z}^{3})\frac{\pa}{\pa z} \right\} \\
  &=&span \left\{ (\phi_{x}^{1} + u \phi_{y}^{1} + v \phi_{z}^{1})\frac{\pa}{\pa x} +  
 (\phi_{x}^{2} + u \phi_{y}^{2} + v \phi_{z}^{2})\frac{\pa}{\pa y} +
  (\phi_{x}^{3} + u \phi_{y}^{3} + v \phi_{z}^{3})\frac{\pa}{\pa z} \right\} \\
  &=&span\left\{ a_{1}\frac{\pa}{\pa x} + a_{2}\frac{\pa}{\pa y} + a_{3}\frac{\pa}{\pa z} \right\}
\end{eqnarray*}
where in the second to last step we divided by 
"$a$" to get that $u = \frac{b}{a}$ and $v = \frac{c}{a}$.
Now, since $\Delta_{1}$ is given by

\

$dy - udx = 0$
\

$dz - vdx = 0$

and since $[dx: dy: dz] = [1: \frac{dy}{dx}: \frac{dz}{dx}]$ we have
for $\Phi \in G_{1}$ that it is given locally as $\Phi^{1}(x,y,z,u,v) = 
(\phi^{1}, \phi^{2}, \phi^{3}, \tilde{u}, \tilde{v})$ where

\

$\tilde{u} = \frac{a_{2}}{a_{1}} = \frac{\phi_{x}^{2} + u \phi_{y}^{2} + v \phi_{z}^{2}}{\phi_{x}^{1} + u \phi_{y}^{1} + v \phi_{z}^{1}}$

\

$\tilde{v} = \frac{a_{3}}{a_{1}} = \frac{\phi_{x}^{3} + u \phi_{y}^{3} + v \phi_{z}^{3}}{\phi_{x}^{1} + u \phi_{y}^{1} + v \phi_{z}^{1}}$

{\it Level $2$:}
At level $2$ we are looking at the class $RV$ which consists of a single orbit.  This means that we can pick any point in the class
$RV$ as our representative. We will pick our point to be $p_{2} = (p_{1}, \ell_{1})$ with $\ell_{1} \subset \Delta_{1}(p_{1})$ to be the vertical line 
$\ell_{1} =[dx: du: dv] = [0: 1: 0]$.  Now, we will let $G_{2}$ be the set of symmetries from $G_{1}$ that fix the vertical line $\ell_{1} = [0: 1: 0]$ in $\Delta_{1}(p_{1})$, meaning we want $\Phi^{1}_{\ast}([0: 1: 0]) = [0: 1: 0]$ for all $\Phi \in G_{2}$.
Then this says that $\Phi^{1}_{\ast}([dx_{| \ell_{1}}: du_{| \ell_{1}}: dv_{| \ell_{1}}]) = \Phi^{1}_{\ast}([0: 1: 0]) = 
[0: 1: 0] = [d \phi^{1}_{| \ell_{1}}: d \tilde{u}_{| \ell_{1}}: d \tilde{v}_{| \ell_{1}}]$.  When we fix this direction it might
yield some new information about the component functions for the $\Phi$ in $G_{2}$.

\

$\bullet$ $d \phi^{1}_{| \ell_{1}} = 0$.
\

$d \phi^{1} = \phi_{x}^{1} dx + \phi_{y}^{1} dy + \phi_{z}^{1} dz$ and that, based on the above, 
$d \phi^{1}\rest{\ell_{1}} = 0$ and can see that we will not gain any new information about the component functions
for $\Phi \in G_{2}$.

\

$\bullet$ $d \tilde{v}_{| \ell_{1}} = 0$
\ 

$d \tilde{v} = d(\frac{a_{3}}{a_{1}}) = \frac{da_{3}}{a_{1}} - \frac{(da_{1})a_{3}}{a_{1}^{2}}$
\

First notice that when we evaluate at $(x,y,z,u,v) = (0,0,0,0,0)$ that $a_{3} = 0$ and since we are setting 
$d \tilde{v}\rest{\ell_{1}} = 0$ then $d a_{3} \rest{\ell_{1}}$ must be equal to zero.  We calculate that
$$da_{3} = \phi_{xx}^{3}dx + \phi_{xy}^{3}dy + \phi_{xz}^{3}dz + \phi_{y}^{3}du + u (d \phi_{y}^{3}) +
          \phi_{z}^{3}dv + v(d \phi_{z}^{3})$$
then when we evaluate we get
$$d a_{3} \rest{\ell_{1}} = \phi^{3}_{y}(\z)du \rest{\ell_{1}} = 0 $$
since $du \rest{\ell_{1}} \neq 0$.  This forces $\phi^{3}_{y}(\z) = 0$. 
\

This gives us the updated Taylor triangle for $\phi^{3}$:

\

\begin{center}
\begin{tabular}{rcccccccccc}
$n=0$:&    &    &    &    &  \xcancel{(0,0,0)} \\\noalign{\smallskip\smallskip}
$n=1$:&    &    &    &  $\xcancel{(1,0,0)}$ & \xcancel{(0,1,0)} & (0,0,1)\\\noalign{\smallskip\smallskip}
$n=2$:&    &    &  (2,0,0) & (1,1,0)   &  (1,0,1) & (0,1,1)    &  (0,0,2)\\\noalign{\smallskip\smallskip}
\end{tabular}
\end{center} 

\

We have determined some of the properties about elements in $G_{2}$ and now we will see what these elements look like locally.
We look at $\Phi^{1}_{\ast}(\ell)$ for $\ell \subset \Delta_{1}$, near the vertical hyperplane in $\Delta_{1}$, and is of the 
form $\ell = a Z^{1} + b \frac{\partial}{\partial u} + 
c \frac{\partial}{\partial v}$ with $a,b,c \in \R$ and $b \neq 0$ with $Z^{1} = u \frac{\partial}{\partial y} + v \frac{\partial}{\partial z} + \frac{\partial}{\partial x}$.  This gives that
\

\[
\Phi^{1}_{\ast}(\ell) =
\begin{pmatrix}
 \phi_{x}^{1}                        &  \phi_{y}^{1}  &  \phi_{z}^{1}  &  0  &  0 \\
 \phi_{x}^{2}                        &  \phi_{y}^{2}  &  \phi_{z}^{2}  &  0  &  0 \\
 \phi_{x}^{3}                        &  \phi_{y}^{3}  &  \phi_{z}^{3}  &  0  &  0 \\
 \frac{\partial \tilde{u}}{\partial x} & \frac{\partial \tilde{u}}{\partial y} & \frac{\partial \tilde{u}}{\partial z} &
 \frac{\partial \tilde{u}}{\partial u} & \frac{\partial \tilde{u}}{\partial v} \\
 \frac{\partial \tilde{v}}{\partial x} & \frac{\partial \tilde{v}}{\partial y} & \frac{\partial \tilde{v}}{\partial z} &
 \frac{\partial \tilde{v}}{\partial u} & \frac{\partial \tilde{v}}{\partial v} 
 \end{pmatrix}
\begin{pmatrix}
 a \\
 au \\
 av \\
 b \\
 c
\end{pmatrix}
\]
\

\begin{eqnarray*}
  = &  & span \{ (a \phi_{x}^{1} + au \phi_{y}^{1} + av \phi_{z}^{1})\frac{\pa}{\pa x} \\
    & + & (a \frac{\partial \tilde{u}}{\partial x} + a u \frac{\partial \tilde{u}}{\partial y} + a v \frac{\partial \tilde{u}}{\partial z} +
          b \frac{\partial \tilde{u}}{\partial u}  + c \frac{\partial \tilde{u}}{\partial v})\frac{\pa}{\pa u} \\
    & + & (a \frac{\partial \tilde{v}}{\partial x} + a u \frac{\partial \tilde{v}}{\partial y} + a v \frac{\partial \tilde{v}}{\partial z} +
          b \frac{\partial \tilde{v}}{\partial u} + c \frac{\partial \tilde{v}}{\pa v})\frac{\pa}{\pa v} \}
\end{eqnarray*}

\begin{eqnarray*}
  = &  &  span \{ (u_{2} \phi_{x}^{1} + uu_{2} \phi_{y}^{1} + vu_{2} \phi_{z}^{1})\frac{\pa}{\pa x} \\
    & + & ( u_{2} \frac{\partial \tilde{u}}{\partial x} + u u_{2} \frac{\partial \tilde{u}}{\partial y} + 
          v u_{2} \frac{\partial \tilde{u}}{\partial z} + \frac{\partial \tilde{u}}{\partial u} + 
          v_{2} \frac{\partial \tilde{u}}{\partial v})\frac{\pa}{\pa u} \\
    & + & ( u_{2} \frac{\partial \tilde{v}}{\partial x} + u u_{2} \frac{\partial \tilde{v}}{\partial y} +  vu_{2} \frac{\partial \tilde{v}}{\partial z}            + \frac{\partial \tilde{v}}{\partial u} + v_{2} \frac{\partial \tilde{v}}{v})\frac{\pa}{\pa v} \}
\end{eqnarray*}

$= span \{ b_{1}\frac{\pa}{\pa x} +  b_{2}\frac{\pa}{\pa u} + b_{3}\frac{\pa}{\pa v} \}$.  Notice that we have only paid attention to the
$x$, $u$, and $v$ coordinates since $\Delta_{1}^{\ast}$ is framed by $dx$, $du$, and $dv$.
Since $u_{2} = \frac{dx}{du}$ and $v_{2} = \frac{dv}{du}$ we get that

\

$\tilde{u}_{2} = \frac{b_{1}}{b_{2}} = \frac{u_{2} \phi_{x}^{1} + uu_{2} \phi_{y}^{1} + vu_{2} \phi_{z}^{1}}{u_{2} \frac{\partial \tilde{u}}{\partial x} + u u_{2} \frac{\partial \tilde{u}}
{\partial y} + v u_{2} \frac{\partial \tilde{u}}{\partial z} + 
 \frac{\partial \tilde{u}}{\partial u} + v_{2} \frac{\partial \tilde{u}}{\partial v}}$

\

$\tilde{v}_{2} = \frac{b_{3}}{b_{2}} = \frac{u_{2} \frac{\partial \tilde{v}}{\partial x} + u u_{2} \frac{\partial \tilde{v}}{\partial y} +  vu_{2} \frac{\partial \tilde{v}}{\partial z} +
  \frac{\partial \tilde{v}}{\partial u} + v_{2} \frac{\partial \tilde{v}}{\partial v}}
{u_{2} \frac{\partial \tilde{u}}{\partial x} + u u_{2} \frac{\partial \tilde{u}}
{\partial y} + v u_{2} \frac{\partial \tilde{u}}{\partial z} + 
 \frac{\partial \tilde{u}}{\partial u} + v_{2} \frac{\partial \tilde{u}}{\partial v}}$
 
\

This tells us what the new component functions $\tilde{u}_{2}$ and $\tilde{v}_{2}$ are for $\Phi^{2}$.

{\it Level $3$:} 
At level $3$ we are looking at the class $RVV$.  We know from our work on the third level that there will be only one orbit
within this class.  This means that we can pick any point in the class
$RVV$ as our representative.  We will pick the point to be $p_{3} = (p_{2}, \ell_{2})$ with $\ell_{2} \subset \Delta_{2}$ to be the vertical line 
$\ell_{2} =[du: du_{2}: dv_{2}] = [0: 1: 0]$.  Now, we will let $G_{3}$ be the set of symmetries from $G_{2}$ that fix the vertical line $\ell_{2} = [0: 1: 0]$ in 
$\Delta_{2}$, meaning we want $\Phi^{2}_{\ast}([0: 1: 0]) = [0: 1: 0] = 
[ d \tilde{u}\rest{\ell_{3}}: d \tilde{u}_{2}\rest{\ell_{3}}: d \tilde{v}_{2}\rest{\ell_{3}} ]$ for all $\Phi \in G_{3}$.
Since we are taking $du\rest{\ell_{3}} 0$ and $dv_{2}\rest{ \ell_{3}} = 0$, with  
$du_{2}\rest{\ell_{3}} \neq 0$ we need to look at $d \tilde{u}\rest{\ell_{3}} = 0$ and $d \tilde{v}_{2}\rest{\ell_{3}} = 0$
to see if these relations will gives us more information about the component functions of $\Phi$.

\

$\bullet$ $d \tilde{u}\rest{\ell_{3}} = 0$.  
\

$d \tilde{u} = d(\frac{a_{2}}{a_{1}}) = \frac{da_{2}}{a_{1}} - \frac{a_{2} da_{1}}{a_{1}^{2}}$ and can see that 
$a_{2}(p_{2}) = 0$ and that 
\

$da_{2}\rest{\ell _{3}} = \phi_{xx}^{2}dx\rest{\ell_{3}} + \phi_{xy}^{2} dy \rest{ \ell_{3}} + \phi_{xz}^{2} dz\rest{ \ell_{3}}
+ \phi_{y}^{2} du \rest{\ell_{3}} + \phi_{z}^{2} dv\rest{\ell_{3}} = 0$.  Since all of the differentials are going to be equal 
to zero when we put the line $\ell_{3}$ into them that we will not gain any new information about the $\phi^{i}$'s.

\

$\bullet$ $d \tilde{v}_{2}\rest{\ell_{3}} = 0$
\

$d \tilde{v}_{2} = d(\frac{b_{3}}{b_{2}}) = \frac{db_{3}}{b_{2}} - \frac{b_{3} db_{2}}{b_{2}^{2}}$.  When we evaluate
we see that $b_{3}(p_{2}) = 0$ since $\frac{\pa \tilde{v}}{\pa u}(p_{2}) = \phi^{3}_{y}(\z) = 0$, which means 
that we only need to look at $\frac{db_{3}}{b_{2}}$.
\begin{eqnarray*}
db_{3} & = & d(u_{2} \frac{\pa \tilde{v}}{\pa x} + u_{2} u \frac{\pa \tilde{v}}{\pa z} + \frac{\pa \tilde{v}}{\pa u} 
+ v_{2} \frac{\pa \tilde{v}}{\pa v}) \\
        & = & \frac{\pa \tilde{v}}{\pa x} du_{2} + u_{2}(d \frac{\pa \tilde{v}}{\pa x}) + u \frac{\pa \tilde{v}}{\pa y} du_{2} \\ 
        & + & u_{2}\frac{\pa \tilde{v}}{\pa y} du + u_{2} u (d \frac{\pa \tilde{v}}{\pa y}) +
  v \frac{\pa \tilde{v}}{\pa z} du_{2} \\
        & + & u_{2} \frac{\pa \tilde{v}}{\pa z} dv + u_{2}v(d\frac{\pa \tilde{v}}{\pa z}) + 
  \frac{\pa \tilde{v}}{\pa u \pa x} dx \\
        & + & \frac{\pa \tilde{v}}{\pa u \pa y} dy + \frac{\pa \tilde{u}}{\pa u \pa z} dz +
  \frac{\pa \tilde{v}}{\pa v}dv_{2} + v_{2}(d \frac{\pa \tilde{v}}{\pa v}).
\end{eqnarray*}  
\

Then evaluating we get
\

$db_{3} \rest{\ell_{3}} = \frac{\pa \tilde{v}}{\pa x}(p_{3})du_{2} \rest{\ell_{3}} = 0$ and since $du_{2} \rest{\ell_{3}} \neq 0$
it forces $\frac{\pa \tilde{v}}{\pa x}(p_{3}) = 0$.  We have $\frac{\pa \tilde{v}}{\pa x}(p_{3}) = 
\frac{\phi^{3}_{xx}(\z)}{\phi^{1}_{x}(\z)}
- \frac{\phi^{1}_{xx}(\z) \phi^{3}_{x}(\z)}{\phi^{1}_{x}(\z)}$ and that $\phi^{3}_{x}(\z) = 0$ to give that
$\frac{\pa \tilde{v}}{\pa x}(p_{3}) = \frac{\phi^{3}_{xx}(\z)}{\phi^{1}_{x}(\z)} = 0$ which forces 
$\phi^{3}_{xx}(\z) = 0$, which gives us information about $\Phi^{3}$.
This gives us the updated Taylor triangle for $\phi^{3}$:

\begin{center}
\begin{tabular}{rcccccccccc}
$n=0$:&    &    &    &    &  \xcancel{(0,0,0)} \\\noalign{\smallskip\smallskip}
$n=1$:&    &    &    &  $\xcancel{(1,0,0)}$ & \xcancel{(0,1,0)} & (0,0,1)\\\noalign{\smallskip\smallskip}
$n=2$:&    &    &  \xcancel{(2,0,0)} & (1,1,0)   &  (1,0,1) & (0,1,1)    &  (0,0,2)\\\noalign{\smallskip\smallskip}
\end{tabular}
\end{center} 

\

Now, our goal is to look at how $\Phi^{3}_{\ast}$'s act on the distribution $\Delta_{3}(p_{3})$ to determine the number
of orbits within the class $RVVV$.
In order to do so we will need to figure out what the local component functions, call them $\tilde{u}_{3}$ and $\tilde{v}_{3}$, for 
$\Phi^{3}$, where $\Phi \in G_{3}$, will look like. To do this we will again look at $\Phi^{2}_{\ast}$ applied to a line $\ell$ that is
near the vertical hyperplane in $\Delta_{2}$.
\

Set $\ell = aZ^{(2)}_{1} + b \frac{\pa}{\pa u_{2}} + c \frac{\pa}{\pa v_{2}}$ for $a,b,c \in \R$ and $b \neq 0$ where
$Z^{(2)}_{1} = u_{2}(u \frac{\pa}{\pa y} + v \frac{\pa}{\pa z} + \frac{\pa}{\pa x}) + \frac{\pa}{\pa u} + v_{2} \frac{\pa}{\pa v}$.
This gives that 

\

\[
\Phi^{2}_{\ast}(\ell) = 
\begin{pmatrix}
\phi^{1}_{x} & \phi^{1}_{y}  & \phi^{1}_{z} & 0 & 0 & 0 & 0 \\
\phi^{2}_{x} & \phi^{2}_{y} & \phi^{2}_{z} & 0 & 0 & 0 & 0 \\
\phi^{3}_{x} & \phi^{3}_{y} & \phi^{3}_{z} & 0 & 0 & 0 & 0 \\
\frac{\pa \tilde{u}}{\pa x} & \frac{\pa \tilde{u}}{\pa y} & \frac{\pa \tilde{u}}{\pa z} &
\frac{\pa \tilde{u}}{\pa u} & \frac{\pa \tilde{u}}{\pa v} & 0 & 0 \\
\frac{\pa \tilde{v}}{\pa x} & \frac{\pa \tilde{v}}{\pa y} & \frac{\pa \tilde{v}}{\pa z} &
\frac{\pa \tilde{v}}{\pa u} & \frac{\pa \tilde{v}}{\pa v} & 0 & 0 \\
\frac{\pa \tilde{u}_{2}}{\pa x} & \frac{\pa \tilde{u}_{2}}{\pa y} & \frac{\pa \tilde{u}_{2}}{\pa z} &
\frac{\pa \tilde{u}_{2}}{\pa u} & \frac{\pa \tilde{u}_{2}}{\pa v} & 
\frac{\pa \tilde{u}_{2}}{\pa u_{2}} & \frac{\pa \tilde{u}_{2}}{\pa v_{2}} \\
\frac{\pa \tilde{v}_{2}}{\pa x} & \frac{\pa \tilde{v}_{2}}{\pa y} & \frac{\pa \tilde{v}_{2}}{\pa z} &
\frac{\pa \tilde{v}_{2}}{\pa u} & \frac{\pa \tilde{v}_{2}}{\pa v} & 
\frac{\pa \tilde{v}_{2}}{\pa u_{2}} & \frac{\pa \tilde{v}_{2}}{\pa v_{2}} 
\end{pmatrix}
\begin{pmatrix}
 au_{2} \\
 a u u_{2} \\
 a v u_{2} \\
 a \\
 a v_{2} \\
 b \\
 c
\end{pmatrix}
\]
  
\

\begin{eqnarray*}
  & = & span\{ (a u_{2} \frac{\pa \tilde{u}}{\pa x} + a u u_{2} \frac{\pa \tilde{u}}{\pa y} + avu_{2} \frac{\pa \tilde{u}}{\pa z}
        + a \frac{\pa \tilde{u}}{\pa u} + av_{2} \frac{\pa \tilde{u}}{\pa v})\frac{\pa}{\pa u} \\
  & + &  (au_{2} \frac{\pa \tilde{u}_{2}}{\pa x} + auu_{2} \frac{\pa \tilde{u}_{2}}{\pa y} + avu_{2} \frac{\pa \tilde{u}_{2}}{\pa z}
         + a \frac{\pa \tilde{u}_{2}}{\pa u} + av_{2} \frac{\pa \tilde{u}_{2}}{\pa v} + 
        b \frac{\pa \tilde{u}_{2}}{\pa u_{2}} + c \frac{\pa \tilde{u}_{2}}{\pa v_{2}})\frac{\pa}{\pa u_{2}} \\
  & + & (au_{2} \frac{\pa \tilde{v}_{2}}{\pa x} + auu_{2} \frac{\pa \tilde{v}_{2}}{\pa y} + avu_{2} \frac{\pa \tilde{v}_{2}}{\pa z} +
         a \frac{\pa \tilde{v}_{2}}{\pa u} + av_{2} \frac{\pa \tilde{v}_{2}}{\pa v}
         + b \frac{\pa \tilde{v}_{2}}{\pa u_{2}} + c \frac{\pa \tilde{v}_{2}}{\pa v_{2}})\frac{\pa}{\pa v_{2}} \}            
\end{eqnarray*}

\begin{eqnarray*}
  & = & span\{ (u_{3} u_{2} \frac{\pa \tilde{u}}{\pa x} + u_{3} u u_{2} \frac{\pa \tilde{u}}{\pa y} + u_{3}vu_{2} \frac{\pa \tilde{u}}{\pa z} + 
        u_{3} \frac{\pa \tilde{u}}{\pa u} + u_{3}v_{2} \frac{\pa \tilde{u}}{\pa v})\frac{\pa}{\pa u} \\
  & + & (u_{3}u_{2} \frac{\pa \tilde{u}_{2}}{\pa x} + u_{3}uu_{2} \frac{\pa \tilde{u}_{2}}{\pa y} + u_{3}vu_{2} \frac{\pa \tilde{u}_{2}}{\pa z} +
         u_{3} \frac{\pa \tilde{u}_{2}}{\pa u} + u_{3}v_{2} \frac{\pa \tilde{u}_{2}}{\pa v} + 
         \frac{\pa \tilde{u}_{2}}{\pa u_{2}} + v_{3} \frac{\pa \tilde{u}_{2}}{\pa v_{2}})\frac{\pa}{\pa u_{2}} \\
  & + &  u_{3}u_{2} \frac{\pa \tilde{v}_{2}}{\pa x} + u_{3}uu_{2} \frac{\pa \tilde{v}_{2}}{\pa y} + u_{3}vu_{2} \frac{\pa \tilde{v}_{2}}{\pa z} +
         u_{3} \frac{\pa \tilde{v}_{2}}{\pa u} + u_{3}v_{2} \frac{\pa \tilde{v}_{2}}{\pa v} + 
         \frac{\pa \tilde{v}_{2}}{\pa u_{2}} + v_{3} \frac{\pa \tilde{v}_{2}}{\pa v_{2}})\frac{\pa}{\pa v_{2}} \} 
\end{eqnarray*}
   
$= span\{ c_{1}\frac{\pa}{\pa u} + c_{2}\frac{\pa}{\pa u_{2}} + c_{3}\frac{\pa}{\pa v_{2}} \}$, since
our local coordinates are given by $[du: du_{2}: dv_{2}] = [\frac{du}{du_{2}}: 1: \frac{dv_{2}}{du_{2}}] = [u_{3}: 1: v_{3}]$ we find that
\

$\tilde{u}_{3} = \frac{c_{1}}{c_{2}} = \frac{u_{3} u_{2} \frac{\pa \tilde{u}}{\pa x} + u_{3} u u_{2} \frac{\pa \tilde{u}}{\pa y} + u_{3}vu_{2} \frac{\pa \tilde{u}}{\pa z} + 
    u_{3} \frac{\pa \tilde{u}}{\pa u} + u_{3}v_{2} \frac{\pa \tilde{u}}{\pa v}}
    {u_{3}u_{2} \frac{\pa \tilde{u}_{2}}{\pa x} + u_{3}uu_{2} \frac{\pa \tilde{u}_{2}}{\pa y} + u_{3}vu_{2} \frac{\pa \tilde{u}_{2}}{\pa z} +
    u_{3} \frac{\pa \tilde{u}_{2}}{\pa u} + u_{3}v_{2} \frac{\pa \tilde{u}_{2}}{\pa v} + 
     \frac{\pa \tilde{u}_{2}}{\pa u_{2}} + v_{3} \frac{\pa \tilde{u}_{2}}{\pa v_{2}}} $

\

$\tilde{v}_{3} = \frac{c_{3}}{c_{2}} = \frac{u_{3}u_{2} \frac{\pa \tilde{v}_{2}}{\pa x} + u_{3}uu_{2} \frac{\pa \tilde{v}_{2}}{\pa y} + u_{3}vu_{2}  \ 
                \frac{\pa \tilde{v}_{2}}{\pa z} +
                 u_{3} \frac{\pa \tilde{v}_{2}}{\pa u} + u_{3}v_{2} \frac{\pa \tilde{v}_{2}}{\pa v} + 
                 \frac{\pa \tilde{v}_{2}}{\pa u_{2}} + v_{3} \frac{\pa \tilde{v}_{2}}{\pa v_{2}}}
                 {u_{3}u_{2} \frac{\pa \tilde{u}_{2}}{\pa x} +   
                 u_{3}uu_{2} \frac{\pa \tilde{u}_{2}}{\pa y} + u_{3}vu_{2} \frac{\pa  
                 \tilde{u}_{2}}{\pa z} +
                 u_{3} \frac{\pa \tilde{u}_{2}}{\pa u} + u_{3}v_{2} \frac{\pa \tilde{u}_{2}}{\pa v} + 
                 \frac{\pa \tilde{u}_{2}}{\pa u_{2}} + v_{3} \frac{\pa \tilde{u}_{2}}{\pa v_{2}}}$.
                 
{\it Level $4$:}
Now that we know what the component functions are for $\Phi^{3}$, with $\Phi \in G_{3}$, we are ready to apply its pushforward
to the distribution $\Delta_{3}$ at $p_{3}$ and figure out how many orbits there are for the class $RVVV$.
We let $\ell = b \frac{\pa}{\pa u_{3}} + c \frac{\pa}{\pa v_{3}}$, with $b,c \in \R$, be a vector in the vertical hyperplane
of $\Delta_{3}(p_{3})$ and we see that 
$$\Phi^{3}_{\ast}(\ell) = 
span \{ (b \frac{\pa \tilde{u}_{3}}{\pa u_{3}}(p_{3}) + c \frac{\pa \tilde{u}_{3}}{\pa v_{3}}(p_{3}))\frac{\pa}{\pa u_{3}} + 
        (b \frac{\pa \tilde{v}_{3}}{\pa u_{3}}(p_{3}) + c \frac{\pa \tilde{v}_{3}}{\pa v_{3}}(p_{3}))\frac{\pa}{\pa v_{3}} \} .$$
This means that we need to see what $\frac{\pa \tilde{u}_{3}}{\pa u_{3}}$, $\frac{\pa \tilde{u}_{3}}{\pa v_{3}}$, 
$\frac{\pa \tilde{v}_{3}}{\pa u_{3}}$, $\frac{\pa \tilde{v}_{3}}{\pa v_{3}}$ are when we evaluate at 
$p_{3} = (x,y,z,u,v,u_{2},v_{2},u_{3},v_{3}) = (0,0,0,0,0,0,0,0,0)$.  This will amount to a somewhat long process, so
we will just first state what the above terms evaluate to and leave the computations for the appendix.
After evaluating we will see that

$\Phi^{3}_{\ast}(\ell) = span \{ (b \frac{(\phi^{2}_{y}(\z))^{2}}{(\phi^{1}_{x}(\z))^{3}})\frac{\pa}{\pa u_{3}} + 
(c \frac{\phi^{3}_{z}(\z)}{(\phi^{1}_{x}(\z))^{2}})\frac{\pa}{\pa v_{3}} \}$, this means
that for $\ell = \frac{\pa}{\pa u_{3}}$ we get $\Phi^{3}_{\ast}(\ell) = 
span \{\frac{(\phi^{2}_{y}(\z))^{2}}{(\phi^{1}_{x}(\z))^{3}}\frac{\pa}{\pa u_{3}} \}$ to 
give one orbit. Then, for $\ell = \frac{\pa}{\pa u_{3}} + \frac{\pa}{\pa v_{3}}$ we see that 
$\Phi^{3}_{\ast}(\ell) = span \{(\frac{(\phi^{2}_{y}(\z))^{2}}{(\phi^{1}_{x}(\z))^{3}})\frac{\pa}{\pa u_{3}} + 
(\frac{\phi^{3}_{z}(\z)}{(\phi^{1}_{x}(\z))^{2}})\frac{\pa}{\pa v_{3}} \}$, and notice
that $\phi^{1}_{x}(\z) \neq 0$, $\phi^{2}_{y}(\z) \neq 0$, and $\phi^{3}_{z}(\z) \neq 0$, but
we can choose them to be anything else to get any vector of the form $b' \frac{\pa}{\pa u_{3}} + c' \frac{\pa}{\pa v_{3}}$ with
$b', c' \neq 0$ to give another, seperate, orbit(Recall that in order for $\ell$ to be a vertical direction, in this case, that it
must be of the form $\ell = b \frac{\pa}{\pa u_{3}} + c \frac{\pa}{\pa v_{3}}$ with $b \neq 0$.)
This means that there is a total of $2$ orbits for the class $RVVV$, as seen in Figure $3$.

\begin{figure}
  \centering
  {\label{fig:orbits}\includegraphics[width = .44 \textwidth]{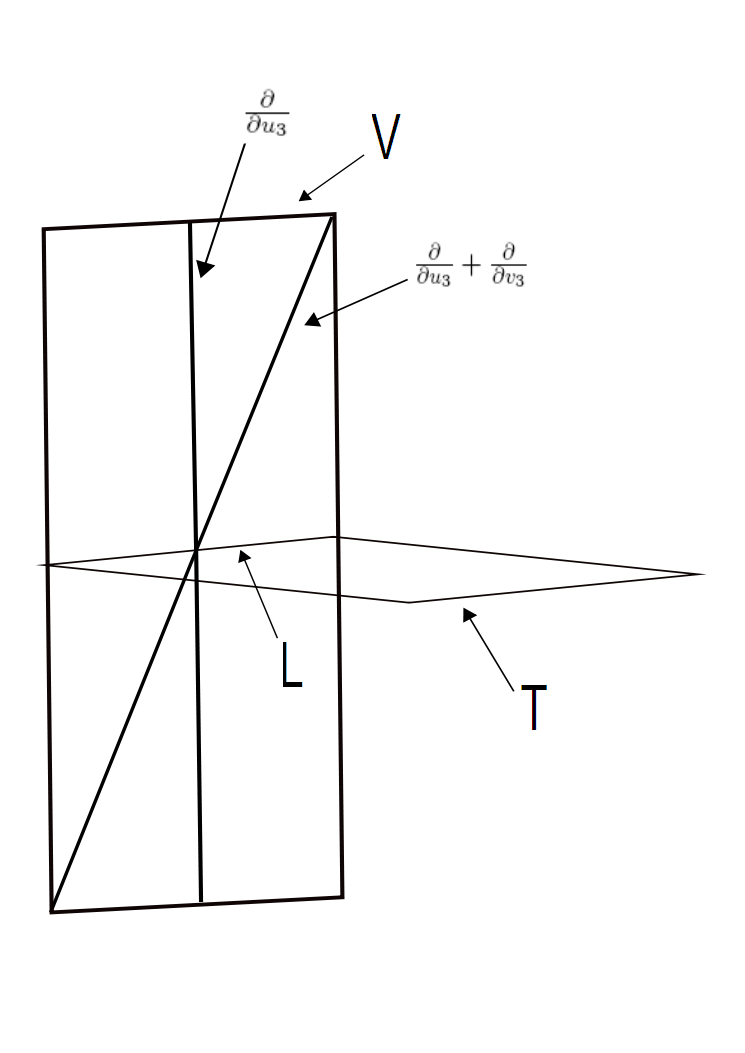}}              
  \caption{Orbits within the class $RVVV$.}
\end{figure}

The classification of the other $RVT$ classes at level $4$ are done in a very similar manner.  The details of these other
calculations will be given in a subsequent work by one of the authors.

\

% ----------------------------------------------------------------
\section{Conclusion} 

We have exbited above a canonical procedure for lifting the action of $Diff(3)$ to the fibered manifold 
$\mathcal{P}^k(2)$, and by a mix of singularity theory of space curves and the representation theory of diffeomorphism groups we were able to completely classify orbits of this extended action for small values of $k$. A cursory glance at our computational methods will convince the reader that these results can nominally be extended to higher values of $k$, say $k\geq 5$ but given the exponential increase in computational effort the direct approach can become somewhat unwieldy.  Progress has been made though to try and extend the classification results of the present paper and we hope to release these findings sometime in the near future.  
In (\cite{castro}) we already called the attention upon a lack of discrete invariants to assist with the classification problem, we hope to return to this problem in future publications. 
\

It came to our attention recently that Respondek and Li \cite{respondek} constructed a mechanical system consisting of $k$-rigid bars moving in $\mathbb{R}^{n+1}$ subjected to a nonholonomic constraint which is equivalent to the Cartan distribution of $J^k(\mathbb{R},\mathbb{R}^{n})$ at regular configuration points. We conjecture that the singular configurations of the $k$-bar will be related to singular Goursat multiflags similar to those presented here, though in Respondek and Li's case the configuration manifold is a tower of $S^n$ fibrations instead of our $\mathbb{P}^n$ tower.

Another research venue, and that to our knowledge has been little explored, is that of understanding how these results could be applied to the geometric theory of differential equations. Let us remind the reader that the spaces $\mathcal{P}^k(2)$, or more generally 
$\mathcal{P}^k(n)$ are vertical compactifications of the jet spaces $J^{k}(\mathbb{R},\mathbb{R}^2)$ and $J^k(\mathbb{R},\mathbb{R}^n)$ respectively.
Kumpera and collaborators (\cite{kumpera2}) have used the geometric theory to study the problem underdetermined systems of systems of ordinary differential equations, but it remains to be explored how our singular orbits can be used to make qualitative statements about the behavior of solutions to singular differential equations.  

%\section{}
%\subsection{}

% -----------------------------------------------------------------
\section{Appendix}

\subsection{A technique to eliminate terms in the short parameterization of a curve germ.}
\

The following technique that we will discuss is outlined in \cite{zariski} on pg. $23$.
Let $C$ be a parameterization of a planar curve germ.  A \textit{short parameterization} of $C$
is of the form

$$C = \begin{cases}
x = t^{n} & \\
y = t^{m} + bt^{\nu_{\rho}} + \Sigma^{q}_{i = \rho + 1} a_{\nu_{i}}t^{\nu_{i}} &\mbox{$b \neq 0$ if $\rho \neq q + 1$}
\end{cases}$$ 

\

\noindent where the $\nu_{i}$, for $i = \rho, \cdots , q$ are positive integers such that $\nu_{\rho} < \cdots < \nu_{q}$
and they do not belong to the semigroup of the curve $C$. Suppose that $\nu_{\rho} + n \in n \Z_{+} + m \Z_{+}$.  
Now, notice that $\nu_{\rho} + n \in m \Z_{+}$ because $\nu_{\rho}$ is not in the semigroup of $C$.  Let
$j \in \Z_{+}$ be such that $\nu_{\rho} + n = (j + 1)m$; notice that $j \geq 1$ since $\nu_{\rho} > m$.  Then set
$a = \frac{bn}{m}$ and \\ $x' = t^{n} + at^{jm} +$( terms of degree $> jm$).  Let 
$\tau^{n} = t^{n} + at^{jm} +$ (terms of degree $> jm$).  From this expression one can show that
$t = \tau - \frac{a}{n} \tau^{jm - n + 1} +$ (terms of degree $> jm - n + 1$), and when we substitute this into the
original expression above for $C$ that 

$$C = \begin{cases}
x' = \tau^{n} \\
y = \tau^{m} + \mbox{terms of degree $> \nu_{\rho}$} 
\end{cases}$$

\

We can now apply semigroup arguments to the above expression for $C$ and see that $C$ has the parametrization

$$C = \begin{cases}
x' = \tau^{n} & \\
y' = \tau^{m} + \Sigma^{q}_{i = \rho + 1} a'_{\nu_{i}} \tau^{\nu_{i}}
\end{cases}$$

\

We can apply the above technique to the two curves $(t^{3}, t^{5} + t^{7}, 0)$ and
$(t^{3}, t^{5} - t^{7}, 0)$ to get that they are equivalent to the curve $(t^{3}, t^{5}, 0)$.

\subsection{Computations for the class $RVVV$.}
\

We will now provide the details that show what the functions $\frac{\pa \tilde{u}_{3}}{\pa u_{3}}$, $\frac{\pa \tilde{u}_{3}}{\pa v_{3}}$, 
$\frac{\pa \tilde{v}_{3}}{\pa u_{3}}$, $\frac{\pa \tilde{v}_{3}}{\pa v_{3}}$ are when we evaluate at 
$p_{4} = (x,y,z,u,v,u_{2},v_{2},u_{3},v_{3}) = (0,0,0,0,0,0,0,0,0)$. 

$\bullet$ $\frac{\pa \tilde{u}_{3}}{\pa u_{3}}$
\

Recall that $\tilde{u}_{3} = \frac{c_{1}}{c_{2}}$.  Then
$\frac{\pa \tilde{u}_{3}}{\pa u_{3}} = \frac{u_{2} \frac{\pa \tilde{u}_{x}}{\pa x} + uu_{2} \frac{\pa \tilde{u}}{\pa y} +
       \frac{\pa \tilde{u}}{\pa u} + v_{2} \frac{\pa \tilde{u}}{\pa}}{c_{2}} - \frac{\frac{\pa c_{2}}{\pa u_{3}} c_{1}}{c^{2}_{2}}$ and
\

$\frac{\pa \tilde{u}_{3}}{\pa u_{3}}(p_{4}) = 
\frac{ \frac{\pa \tilde{u}}{\pa u}(p_{4}) }{ \frac{\pa \tilde{u}_{2}}{\pa u_{2}}(p_{4})}$, since
$c_{1}(p_{4}) = 0$.  We recall that $\frac{\pa \tilde{u}}{\pa u}(p_{4}) = \frac{\phi^{2}_{y}(\z)}{\phi^{1}_{x}(\z)}$, 
$\frac{\pa \tilde{u}_{2}}{\pa u_{2}}(p_{4}) = \frac{\phi^{1}_{x}(\z)}{\frac{\pa \tilde{u}}{\pa u}(p_{4})}$ to give that
\

$\frac{\pa \tilde{u}_{3}}{\pa u_{3}}(p_{4}) = \frac{(\phi^{2}_{y}(\z))^{2}}{(\phi^{1}_{x}(\z))^{3}}$.

\

$\bullet$ $\frac{\pa \tilde{u}_{3}}{\pa v_{3}}$.
\

Have that $\tilde{u}_{3} = \frac{c_{1}}{c_{2}}$, then

\

$\frac{\pa \tilde{u}_{3}}{\pa v_{3}}(p_{4}) = \frac{\frac{\pa c_{1}}{\pa v_{3}}(p_{4})}{c_{2}(p_{4})} - 
\frac{\frac{\pa c_{2}}{\pa v_{3}}(p_{4}) c_{1}(p_{4})}{c^{2}_{2}(p_{4})} =
 0$, since $c_{1}$ is not a function of $u_{3}$ and $c_{1}(p_{4}) = 0$.

\

$\bullet$ $\frac{\pa \tilde{v}_{3}}{\pa u_{3}}$
\

Have that $\tilde{v}_{3} = \frac{c_{3}}{c_{2}}$, then

\

$\frac{\pa \tilde{v}_{3}}{\pa u_{3}} = \frac{u_{2}\frac{\pa \tilde{v}_{2}}{\pa x} + ... + \frac{\pa \tilde{v}_{2}}{\pa u} + 
                                       v_{2} \frac{\pa \tilde{v}_{2}}{\pa v}}{c_{2}} 
                                       - \frac{(u_{2} \frac{\pa \tilde{u}_{2}}{\pa x} + ... + \frac{\pa \tilde{u}_{2}}{\pa u} +
                                       ... + v_{2} \frac{\pa \tilde{u}_{2}}{\pa v})c_{1}}{c^{2}_{2}}$
                                       
\

$\frac{\pa \tilde{v}_{3}}{\pa u_{3}}(p_{4}) = \frac{ \frac{\pa \tilde{v}_{2}}{\pa u} (p_{4})}{ \frac{\pa \tilde{u}_{2}}{\pa u_{2}}(p_{4}) }
                                                      - \frac{ \frac{\pa \tilde{u}_{2}}{\pa u}(p_{4}) \frac{\pa \tilde{v}_{2}}{\pa u_{2}}(p_{4}) }
                                                      { (\frac{\pa \tilde{u}_{2}}{\pa u_{2}}(p_{4}))^{2} }$

\

This means that we will need to figure out what $\frac{\pa \tilde{u}_{2}}{\pa u_{2}}$, $\frac{\pa \tilde{v}_{2}}{\pa u_{2}}$,
and $\frac{\pa \tilde{v}_{2}}{\pa u}$ are
when we evaluate at $p_{4}$.

\

$\circ$ $\frac{\pa \tilde{u}_{2}}{\pa u_{2}}$
\

Can recall from work at the level below that $\frac{\pa \tilde{u}_{2}}{\pa u_{2}}(p_{4}) = 
\frac{\phi^{1}_{x}(\z)}{\frac{\pa \tilde{u}}{\pa u}(p_{4})} = \frac{\phi^{2}_{y}(\z)}{(\phi^{1}_{x}(\z))^{2}}$ since 
$\frac{\pa \tilde{u}}{\pa u}(p_{4}) = \frac{\phi^{2}_{y}(\z) }{\phi^{1}_{x}(\z)}$.

\

$\circ$ $\frac{\pa \tilde{v}_{2}}{\pa u_{2}}$
\

Can recall from work at level $3$ that $\frac{\pa \tilde{v}_{2}}{\pa u_{2}}(p_{4}) = 
\frac{\frac{\pa \tilde{v}}{\pa x}(p_{4})}{\frac{\pa \tilde{u}}{\pa u}(p_{4})} = 0$ since 
$\frac{\pa \tilde{v}}{\pa x}(p_{4}) = \frac{\phi^{3}_{xx}(\z)}{\phi^{1}_{x}(\z)}$
and have that $\phi^{1}_{xx}(\z) = 0$ to give $\frac{\pa \tilde{v}_{2}}{\pa u_{2}}(p_{4}) = 0$.

This gives the reduced expression $\frac{\pa \tilde{v}_{3}}{\pa u_{3}}(p_{4}) = 
\frac{\frac{\pa \tilde{v}_{2}}{\pa u}(p_{4})}{\frac{\pa \tilde{u}_{2}}{\pa u_{2}}(p_{4})}$.
\

$\circ$ $\frac{\pa \tilde{v}_{2}}{\pa u}$
\

Recall that $\frac{\pa \tilde{v}_{2}}{\pa u} = \frac{b_{3}}{b_{2}}$, then we get that

\

$\frac{\pa \tilde{v}_{2}}{\pa u} = 
\frac{ u_{2} \frac{\pa \tilde{v}}{\pa x \pa u} + u_{2} \frac{\pa \tilde{v}}{\pa y} + ... + \frac{\pa \tilde{v}}{\pa^{2} u} + 
v_{2} \frac{\pa \tilde{v}}{\pa v \pa u}}{b_{2}} - \frac{(u_{2} \frac{\pa \tilde{u}}{\pa x \pa u} + ... + \frac{\pa \tilde{u}}{\pa^{2} u} + 
v_{2} \frac{\pa \tilde{u}}{\pa v \pa u})b_{3}}{b_{2}^{2}}$

\

$\frac{\pa \tilde{v}_{2}}{\pa u}(p_{4}) = \frac{\frac{\pa \tilde{v}}{\pa^{2} u}(p_{4})}{ \frac{\pa \tilde{u}}{\pa u}(p_{4})} - 
\frac{\frac{\pa \tilde{u}}{\pa^{2} u}(p_{4})  \frac{\pa \tilde{v})}{\pa u}(p_{4})}{ (\frac{\pa \tilde{u}}{\pa u}(p_{4}))^{2}}$ since 
$b_{2}(p_{4}) = \frac{\pa \tilde{u}}{\pa u}(p_{4})$ and $b_{3}(p_{4}) = \frac{\pa \tilde{v}}{\pa u}(p_{4})$.  In order to 
figure out $\frac{\pa \tilde{v}_{2}}{\pa u}(p_{4})$ will be we need to look at 
$\frac{\pa \tilde{v}}{\pa u}(p_{4})$, $\frac{\pa \tilde{v}}{\pa^{2} u}(p_{4})$, and $\frac{\pa \tilde{u}}{\pa^{2} u}(p_{4})$.

\

$\circ$ $\frac{\pa \tilde{v}}{\pa u}$
\

Recall that $\tilde{v} = \frac{a_{3}}{a_{1}}$ and that 
$\frac{\pa \tilde{v}}{\pa u} = \frac{\phi^{3}_{y}}{a_{1}} - \frac{\phi^{1}_{y} a_{3}}{a^{2}_{1}}$, then

\

$\frac{\pa \tilde{v}}{\pa u}(p_{4}) = \frac{\phi^{3}_{y}(\z)}{\phi^{1}_{x}(\z)} - 
\frac{ \phi^{1}_{y}(\z) \phi^{3}_{x}(\z)}{(\phi^{1}_{x}(\z))^{2}} = 0$ since
$\phi^{3}_{y}(\z) = 0$ and $\phi^{3}_{x}(\z) = 0$.
\

$\circ$ $\frac{\pa \tilde{v}}{\pa^{2} u}$
\

From the above we have $\frac{\pa \tilde{v}}{\pa u} = \frac{\phi^{3}_{y}}{a_{1}} - \frac{\phi^{1}_{y} a_{3}}{a^{2}_{1}}$, then
\

$\frac{\pa \tilde{v}}{\pa^{2} u}(p_{4}) = \frac{0}{a_{1}(p_{4})} - \frac{\phi^{3}_{y}(\z) \phi^{1}_{y}(\z)}{a^{2}_{1}(p_{4})} - 
\frac{\phi^{1}_{y}(\z) \phi^{3}_{y}(\z)}{a^{2}_{1}(p_{4})} + \frac{(\phi^{1}_{y}(\z))^{2} \phi^{3}_{x}(\z)}{a^{3}_{1}(p_{4})} = 0$ since 
$\phi^{3}_{y}(\z) = 0$ and 
$\phi^{3}_{x}(\z) = 0$.
\

We see that we do not need to determine what $\frac{\pa \tilde{u}}{\pa^{2} u}(p_{4})$ is, since 
$\frac{\pa \tilde{v}}{\pa u}$ and $\frac{\pa \tilde{v}}{\pa^{2} u}$ will be zero at $p_{4}$ and give us that
$\frac{\pa \tilde{v}_{3}}{\pa u_{3}}(p_{4}) = 0$.

\

$\bullet$ $\frac{\pa \tilde{v}_{3}}{\pa v_{3}}$.
\

Recall that $\tilde{v}_{3} = \frac{c_{3}}{c_{2}}$, then

\

$\frac{\pa \tilde{v}_{3}}{\pa v_{3}} = \frac{ u_{3}u_{2} \frac{\pa \tilde{v}_{2}}{\pa x \pa v_{3}} + ... + \frac{\pa \tilde{v}_{2}}{\pa v_{2}}  
        }{c_{2}} - \frac{ (u_{3}u_{2} \frac{\pa \tilde{u}_{2}}{\pa x \pa v_{3}} + ... + \frac{\pa \tilde{u}}{\pa v_{2}})c_{3}}{c^{2}_{2}}$
        
\

$\frac{\pa \tilde{v}_{3}}{\pa v_{3}}(p_{4}) = \frac{ \frac{\pa \tilde{v}_{2}}{\pa v_{2}}(p_{4})}{ \frac{\pa \tilde{u}_{2}}{\pa u_{2}}(p_{4})} - 
\frac{ \frac{\pa \tilde{u}}{\pa v_{2}(p_{4})} \frac{\pa \tilde{v}_{2}}{\pa u_{2}}(p_{4})}{(\frac{\pa \tilde{u}_{2}}{\pa u_{2}}(p_{4}))^{2}}$.  This means that we need to look at $\frac{\pa \tilde{v}_{2}}{\pa v_{2}}$, $\frac{\pa \tilde{u}}{\pa v_{2}}$, $\frac{\pa \tilde{v}_{2}}{\pa u_{2}}$, and
$\frac{\pa \tilde{u}_{2}}{\pa u_{2}}$ evaluated at $p_{4}$.

\

$\circ$ $\frac{\pa \tilde{v}_{2}}{\pa v_{2}}$.
\

We recall from an earlier calculation that $\frac{\pa \tilde{v}_{2}}{\pa v_{2}}(p_{4}) = 
\frac{ \frac{\pa \tilde{v}}{\pa v}(p_{4})}{\frac{\pa \tilde{u}}{\pa u}(p_{4})} = \frac{\phi^{3}_{z}(\z)}{\phi^{2}_{y}(\z)}$.

\

$\circ$ $\frac{\pa \tilde{u}}{\pa v_{2}}$.
\

Recall from an earlier calculation that $\frac{\pa \tilde{u}_{2}}{\pa v_{2}}(p_{4}) = 0$.

\

$\circ$ $\frac{\pa \tilde{u}_{2}}{\pa u_{2}}$.
\

Recall that $\tilde{u}_{2} = \frac{b_{1}}{b_{2}}$ and that

$\frac{\pa \tilde{u}_{2}}{\pa u_{2}} = \frac{ \phi^{1}_{x} + u \phi^{1}_{y} + v \phi^{1}_{z}}{b_{2}} - 
\frac{(\frac{\pa \tilde{u}}{\pa x} + u \frac{\pa \tilde{u}}{\pa y} + v \frac{\pa \tilde{u}}{\pa z})b_{1}}{b^{2}_{2}}$, then

\

$\frac{ \pa \tilde{u}_{2}}{\pa u_{2}}(p_{4}) = \frac{\phi^{1}_{x}(\z)}{\frac{\pa \tilde{u}}{\pa u}(p_{4})} = 
\frac{(\phi^{1}_{x}(\z))^{2}}{\phi^{2}_{y}(p_{4})}$.

With the above in mind we see that $\frac{\pa \tilde{v}_{3}}{\pa v_{3}}(p_{4}) = \frac{\phi^{3}_{z}(\z)}{(\phi^{1}_{x}(\z))^{2}}$.

\

This gives that $\Phi^{3}_{\ast}(\ell) = span \{ (b \frac{\pa \tilde{u}_{3}}{\pa u_{3}}(p_{4}) + 
c \frac{\pa \tilde{u}_{3}}{\pa v_{3}}(p_{4}))\frac{\pa}{\pa u_{3}} + 
(b \frac{\pa \tilde{v}_{3}}{\pa u_{3}}(p_{4}) + c \frac{\pa \tilde{v}_{3}}{\pa v_{3}}(p_{4}))\frac{\pa}{\pa v_{3}} \} = 
span \{ (b \frac{(\phi^{2}_{y}(\z))^{2}}{(\phi^{1}_{x}(\z))^{3}})\frac{\pa}{\pa u_{3}} + 
c \frac{\phi^{3}_{z}(\z)}{(\phi^{1}_{x}(\z))^{2}}\frac{\pa}{\pa v_{3}} \}$.

%%%%%%%%%%%%%%%%%%%%%%%%%%%%%%%%%%%%%%%%%%%%
\medskip
\bibliographystyle{alpha}	%alpha style uses authors last name and year
\bibliography{orbits-wWYATTref}

%%%%%%%%%%%%%%%%%%%%%%%%%%%%%%%%%%%%%%%%%%%%%

\end{document}